\documentclass[12pt,a4paper,twoside, notitlepage]{article}

\usepackage{amsmath}
\usepackage{dsfont}
\usepackage{mathe}
\usepackage{amssymb}
\usepackage{enumerate}

\usepackage{graphicx}

\usepackage{fancyhdr}
\usepackage{diss}

\def\hyph{\penalty0\hskip0pt\relax}

\def\slope{{\rm slope}}
\def\Fold{\mathop{\rm Fold}}
\def\star{\mathop{\rm star}}
\def\supp{\mathop{\rm supp}}
\def\cotan{\mathop {\rm cotan}\nolimits}
\def\Conv{\mathop {\rm Conv}}
\def\M{\mathcal{M}}
\def\F{\mathcal{F}}
\def\C{\mathcal{C}}
\def\EB{Euclidean building}
\def\HS{Hadamard space}

\def\T{\mathcal{T}}
\def\K{\mathcal{K}}

\def\comment#1{}

\def\Wlo{Without loss of generality}
\def\wlo{without loss of generality}

\begin{document}

\title{Convex rank 1 subsets of Euclidean Buildings (of type $A_{2}$)}
\author{Andreas Balser}

\maketitle
\begin{abstract}
\noindent
For a \EB ~$X$ of type $A_{2}$, 
we 
classify the 0-dimensional subbuildings $A$ of $\partial_{T}X$ that  occur as
the asymptotic boundary of closed convex subsets. 
In particular, we show that triviality of the holonomy of a
triple (of points of $A$) is (essentially) sufficient.
To prove this, we construct new convex subsets as the union of convex sets.
 \footnote{
\emph{keywords:} \EB s, \HS s}
\end{abstract}

\section{Introduction}

%
%

Recently, convex subsets of symmetric spaces were studied in 
\cite{kleinerLeebGroups}:
 Bruce Kleiner and Bernhard Leeb classify
which convex subsets of $\partial_{T}X$ arise as the asymptotic boundary
of a convex subset $C$ of $X$, which is invariant under a group of
isometries of the 
symmetric space $X$ acting cocompactly on~$C$.

They show (via a careful analysis of the Tits boundary
$\partial_{T}C\subset \partial_{T}X$) that such a set 
$C$ is a symmetric subspace, a subset of a 
rank~1-symmetric subspace, or a product of sets of this kind.

We omit the group action and ask:
Which (closed) $\pi$-convex subsets of $\partial_{T}X$ can
occur as the asymptotic boundary of a convex subset of $X$?

\medbreak

If $\partial_{T}X$ is a spherical building, then
 $\pi$-convex subsets of $\partial_{T}X$ of dimension
at most two have a center or are subbuildings (\cite{centers}). 


In this article,
we restrict our attention to
0-dimensional subbuildings. I.e.~we consider subsets $A\subset
\partial_{T}X$ such that every pair $\eta ,\xi \in A$ satisfies
$\angle_{Tits}(\eta ,\xi )\geq \pi$. 

We call a subset $C\subset X$ of a \HS ~$X$ a 
\emph{convex rank 1}-subset if it is closed, convex, and its
asymptotic boundary  
$\partial_{T}C$ is a 0-dimensional subbuilding
of~$\partial_{T}X$
(see 
also definition \ref{convexRank1Definition}).


We will find that it is necessary for each triple of points
$\eta_{1},\eta_{2},\eta_{3}$ of $A$ to
correspond to an ideal triangle: That is, it is necessary that there
are lines $l_{i,j}\subset X, i,j\in \{1,2,3 \}$ joining the boundary
points $\eta_{i},\eta_{j}$, which are pairwise strongly
asymptotic (at the common endpoint). 
Formally speaking, it is necessary that the \emph{holonomy map} of
this triple has a fixed point.
For more details on holonomy, see Section \ref{holonomySection}.

In the case where $X$ is a \EB~of type $A_{2}$, we find that this
holonomy condition is (essentially) sufficient; this is shown by
constructing new convex sets as the union of convex subsets of $X$.
Our main theorem is:

\begin{theorem}\label{convBuild}
Let $X$ be a \EB ~of type $A_{2}$, and let $A\subset
\partial_{T}X$ be a finite 0-dimensional subbuilding of its boundary.
Then there exists a convex rank 1-subset $C\subset X$ such that
$\partial_{T}C \supset A$ if and only if each triple of points of~$A$
corresponds to an ideal triangle.

If  $A$ is infinite, the claim holds under an
additional necessary assumption ($A$ needs to be ``good'', see Definition
\ref{goodDef}). 
\end{theorem}

In the proof, we construct a convex set $C$ which satisfies
$\partial_{T}C=\bar{A}$, 
where $\bar{A}$ is the closure of $A\subset \partial_{\infty X}$, the
asymptotic boundary of $X$ with the cone topology.

The idea to examine this question, and the technique of using
holonomy, are due to Bruce Kleiner and Bernhard Leeb, who can classify
the possible boundaries of convex rank 1-subsets of 
$\R H^2 \times \R H^2$.

This topic is also related to the work \cite{schroederNegConvHulls} of
Hummel, Lang and Schroeder:  
They show that in a CAT(-1)-space, the convex hull of finitely many
closed convex sets lies in a (finite) tubular neighbor of this union. 

We examine how to generalize this to
\EB s of type $A_{2}$, where 
the starting blocks may be considered lines, or rather tripods, see
Prop.~\ref{existenceOfTripods} 
 (observe that every subset of the asymptotic boundary
of a CAT(-1)-space is a 0-dimensional subbuilding).

Theorem \ref{convBuild} forms a contrast to the following result: 

\begin{theorem}\label{convSymm}
Let $C$ be a convex rank 1-subset of $M:=SL(3,\R)/SO(3,\R)$. 
Then we have $\partial_{T}C\subset \partial_{T}\H^2$ for
a suitable isometric embedding (up to rescaling) 
of the hyperbolic plane
$\H^2\hookrightarrow M$. 
\end{theorem}

This more expected result is in line with the results of
\cite{kleinerLeebGroups}, and makes it hard to predict Theorem
\ref{convBuild}. A proof of Theorem \ref{convSymm} can be found in
\cite[ch.~IV]{balserDiss}.

\medbreak

For 2-dimensional \EB s, the Coxeter complex $A_{2}$ is special: 

%

In the other three 2-dimensional Coxeter complexes
($B_{2},G_{2}$, and $A_{1}\times A_{1}$), the situation is different:
In contrast to $A_{2}$ (where every holonomy map is orientation
preserving),
the holonomy map of a pair
of antipodal points is orientation reversing in the other cases. 

In these cases, 
existence of tripods is the essential question: 
If one could show that tripods exist (as in
\ref{existenceOfTripods}), then the orientation-reversing property of
holonomy maps leads immediately to the conclusion
that
there exists an isometrically embedded tree.

\tableofcontents

\section{Layout of the Paper}
In section \ref{HSsection}, we introduce basic facts about \HS s and
\EB s needed later on.

In Proposition \ref{convexityLocalProperty}, we show that for a
connected, closed subset $C$ of a Hadamard space, convexity is a local
property. This result feels like a version of the Hadamard-Cartan
theorem; however, it is not an immediate consequence, because we do
not know that $C$ is a geodesic space itself.
%

This proposition will be the tool to show that the sets we construct
in the proof of our main theorem are convex.

In section \ref{busemannSection}, we examine under which conditions
Busemann functions agree, and in which cases unions of horoballs are
(locally) convex. The lemmas in this section are formulated generally
for \EB s, and we hope that they will be useful in the study of
convex rank 1-subsets of higher-dimensional buildings.

Section \ref{geometrySection} contains geometric lemmas about
buildings of type $A_{2}$: We exclude the existence of triangles
$\Delta (a,b,c)$ with
certain properties. To formulate it positively, we show that under
certain circumstances, the starting direction $\overrightarrow{xc}$
(for $x\in \overline{ab}$) always points in ``roughly the same
direction''. 

Starting with section \ref{necessarySection}, we move directly towards
the proof of Theorem~\ref{convBuild}. From there on, $X$ always stands
for a \EB~of type $A_{2}$.

In section \ref{necessarySection}, we examine necessary conditions for
$A\subset \partial_{T}X$ to lie in the boundary of a convex rank
1-set. In particular, we show that for every triple of boundary
points, a tripod has to exist (Proposition
\ref{existenceOfTripods}). We call $A\subset \partial_{T}X$ an S-set
if it satisfies this condition.

If one knows (or expects) Theorem \ref{convSymm}, one might also
expect 
that for a building of type $A_{2}$, every convex rank 1-subset is
essentially a tree. This turns out to be wrong. However, in section
\ref{treeSection},  we obtain a
tree $\T$ as a quotient of a subset of $X$ naturally associated to
the S-set $A \subset \partial_{T}X$.

In section \ref{thickeningSection}, we ``thicken'' tripods; i.e.~we
search for convex rank 1-subsets of $X$ containing a given tripod.
This motivates the definition of the convex set $\K$ in the 
section which follows, and introduces the techniques for proving
convexity. 

The last section \ref{finalSection}, finally, presents the proof of Theorem
\ref{convBuild}: Given a good S-set $A\subset \partial_{T}X$, we
consider the associated tree $\T$.  For every point $[x]\in \T$, we
define a closed  
convex subset $\K_{[x]}$ of $X$, and we show that the (closure of the)
union $\K$ 
of these sets is convex. In a last step, we obtain a subset $\bar{\C}$
of $\bar{\K}$ which can easily be seen to be convex rank 1, and
satisfy $\partial_{T}\bar{\C}=\bar{A}$, where $\bar{A}$ is the closure of
$A\subset \partial_{\infty}X$ in the cone topology.


\subsection*{Acknowledgements}
In this article, I present the main result of my PhD thesis
\cite{balserDiss}: 

Hence, first of all, I would like to express my gratitude to Bernhard Leeb
for supervising me, 
posing the topic and turning me towards the right questions.

It is a pleasure to thank my
discussion-mates Robert Kremser and Carlos Ramos Cuevas.

I would like to thank Alexander Lytchak for fruitful conversations about
\cite{balserPolygons}, and for the collaboration on
\cite{buildingLike} and \cite{centers}. This teamwork gave me a lot of
insights and 
inspiration for the questions discussed here.


\section{\HS s}
\label{HSsection}
We will use the language of non-positively curved metric spaces, as
developed in \cite{ballmann}. 

Throughout, let $X$ be a Hadamard space, unless otherwise stated.
We will use the terms \HS~and CAT(0)-space synonymously; i.e., we
impose completeness on every CAT(0)-space.
Note that a \HS\ need not be locally compact.

Recall that $X$ has a \emph{boundary at infinity}
$\partial_{\infty}X$, which is given by equivalence classes of rays,
where two (unit-speed) rays are equivalent if their distance is
bounded. 

In particular, we will use Busemann functions $b_{\eta}$ associated to
an asymptotic boundary point $\eta \in \partial_{\infty}X$.
A Busemann function measures (relative) distance from a point at
infinity, and is determined up to an additive constant
only. Busemann functions are convex (along any geodesic) and
1-Lipschitz. 

Geodesics, rays, and geodesic segments are always assumed to be
pa\-ra\-me\-trized by unit speed (i.e.\ they are isometric embeddings). 

For a line $l$ in $X$, there is the space $P_{l}$ of parallel
lines. $P_{l}$ splits as a product $P_{l}\cong l\times CS(l)$, where
$CS(l)$ is a \HS \ again. 

For points $x,\xi$ with $x\in X$, $\xi \in X\cup
\partial_{\infty}X=:\bar{X}$, 
and $t\geq 0$ (if $\xi \in X$, let $t\leq d(x,\xi )$), 
we let $\overline{x\xi}(t)$ denote the point on the segment/ray
$\overline{x\xi}$ at distance $t$ from $x$. When we denote a ray by
$\overline{o\eta}$, we order the points such that $o\in X$ and $\eta
\in \partial_{\infty}X$. 

In our notation, $B_{r}(o):=\{x\in X \st d(x,o)\leq r\}$, for $o\in
X,r\geq 0$; i.e.\ balls in \HS s are always assumed closed.

Whenever $C$ is a closed convex subset of a \HS ~$X$, then
$\pi_{C}:X\rightarrow C$ denotes the nearest-point projection (see
\cite[II.2.4]{bridsonHaefliger}). 

\subsection{Angles, spaces of directions, and Tits distance}
Let $o\in X$ be a point in a \HS , and let $\eta ,\xi \in
\partial_{\infty}X$. Let $c,c'$ be the rays
$\overline{o\eta},\overline{o\xi}$. For points $c(t),c'(t')$, one can
consider the \emph{Euclidean comparison triangle} corresponding to the
points $o,c(t),c'(t')$, i.e.\ the Euclidean triangle with side-lengths
$d(o,c(t)), d(c(t),c'(t')), d(c'(t'),o)$ (which is well-defined up to
isometries of the Euclidean plane). The \emph{comparison angle}
between $c(t)$ and $c'(t')$ at $o$ is the angle of the comparison
triangle at the point corresponding to~$o$. It is denoted by
$\tilde{\angle}_{o}(c(t),c'(t'))$. 

We have the following monotonicity property: 
\[
0<t\leq s\text{ and }0<t'\leq s'\text{ implies }
\tilde{\angle}_{o}(c(t),c'(t'))\leq
\tilde{\angle}_{o}(c(s),c'(s')).
\]

From this, one can deduce a notion of angle between geodesic segments
and rays: 
\[
\angle_{o}(\eta,\xi  ):=\lim_{t,t'\rightarrow
0}\tilde{\angle}_{o}(c(t),c'(t'))\in [0,\pi ],
\]
and an ``angle at infinity'', the \emph{Tits angle} between boundary
points
\[
\angle (\eta,\xi  ):=\angle_{Tits}(\eta ,\xi ):=\lim_{t,t'\rightarrow
\infty}\tilde{\angle}_{o}(c(t),c'(t'))	\in [0,\pi ].
\]
It is easy to see that the Tits angle between $\eta ,\xi$ does not
depend on the chosen 
basepoint~$o$. 
The length metric induced on $\partial_{\infty}X$ by $\angle$ is
called \emph{Tits distance} $Td$, and makes $\partial_{\infty}X$ a
CAT(1)-space. If one wants to emphasize that  the Tits
distance and corresponding topology on $\partial_{\infty}X$ is
considered, this 
space is sometimes called $\partial_{T}X$. We will use these
expressions synonymously (and we usually consider the Tits topology).
If the Tits angle (between $\eta ,\xi$) is less than $\pi$, there is a 
unique geodesic $\overline{\eta \xi}\subset \partial_{\infty}X$
connecting them. 

Similarly, the space of directions $\Sigma_{o}(X)$, i.e.\ the completion of
the space 
of starting directions of geodesic segments initiating in
$o$ (modulo the equivalence of directions enclosing a zero angle), can
be regarded as a CAT(1)-space. For $o\in X,x\in \bar{X}$, we let
$\overrightarrow{ox}\in \Sigma_{o}(X)$ be the starting direction of
the segment~$\overline{ox}$.

We call a subset $C\subset B$ of a CAT(1)-space $B$ convex if it is
$\pi$-convex, i.e.\ if all pairs of points of distance less than $\pi$
can be 
joined by a geodesic.

\subsection{Strong asymptote classes and holonomy}
\label{holonomySection}
Two rays $\overline{o\eta},\overline{x\eta}$ in a \HS~$X$ are called
strongly asymptotic if 
\[
d_{\eta}(\overline{o\eta},\overline{x\eta}):=\lim_{t\rightarrow
\infty}d(\overline{o\eta}(t),\overline{x\eta})=0.
\]
This defines an
equivalence relation on the set of rays asymptotic to $\eta$. 
The metric completion $X_{\eta}$ of this set of equivalence classes is called 
the space of strong asymptote classes at $\eta$. It is a \HS~again
(see \cite{karpelevic}, \cite[sect.\ 2.1.3]{leebHabil}).

Now assume that $X$ is a symmetric space or a \EB , and consider two
antipodal points $\eta ,\xi \in \partial_{T}X$.
It is well known that the parallel set of $(\eta ,\xi )$, i.e.~all the
lines with asymptotic endpoints $\eta ,\xi$,
 represents all the strong asymptote classes at $\eta$ and
at $\xi$.

This induces a natural isometry $h_{\eta ,\xi}:X_{\eta}\rightarrow
X_{\xi}$. 

%

Such a map (and composition of such maps) is called a \emph{holonomy
map} (see \cite[ch.~3]{leebHabil}).

\subsection{\EB s}
We will also need some Euclidean building geometry. For an
introduction, we refer to \cite[sect.~4]{kleinerLeeb}.
A brief introduction of the notation we use can be found in
\cite[sect.~2.4]{kapovichLeebMillson}.
Note that in particular, a \EB \ is a Hadamard space.

A 1-dimensional \EB~is called a \emph{tree}.

In a \EB, we call a geodesic segment regular, if all its interior
points are regular.

The boundary at infinity of a \EB \ $X$ of rank $n$ is a spherical
building of dimension $n-1$; 
we refer to \cite[sect.~3]{kleinerLeeb} for an introduction.

We will use that a spherical building $B$ is a spherical simplicial
complex, where all the simplices are isometric to a spherical polytope
$\Delta$ (in particular, $\Delta$ tesselates $S^{n-1}$), which is the
\emph{spherical Weyl chamber} of the building. Apartments
(i.e.\ isometrically embedded copies $S^{n-1}$) intersect in (unions
of) Weyl
chambers.
There is a natural map $B\rightarrow \Delta$, and the
image of a point is called its \emph{type}.

\section{Convexity is a local property in CAT(0)-spaces}
\label{convexitySection}
Let $\varepsilon >0$. A subset $C$ of a CAT(0)-space $X$ is called
$\varepsilon$-locally convex, if for all $x\in C$, the set
$B_{\varepsilon}(x)\cap C$ is convex. Note that a convex set is
$\varepsilon$-locally convex for all $\varepsilon >0$. 
Since an $\varepsilon$-locally convex set is locally path-connected,
path-connectedness and connectedness are equivalent for
$\varepsilon$-locally convex sets.

We show that if $C$ is closed and connected, then one
$\varepsilon$ suffices to make sure that $C$ is convex:

\begin{proposition}\label{convexityLocalProperty}
Let $\varepsilon>0$, and $X$ be a CAT(0)-space. Let $C\subset X$ be a
connected, closed, 
$\varepsilon$-locally convex set. Then $C$ is convex.
\end{proposition}

Observe that this claim is similar in nature to the Hadamard-Cartan
theorem, saying that a geodesic space which is simply connected and
locally CAT(0), is actually globally CAT(0). In our case, we do not
know whether~$C$ is a geodesic space, so this proposition is not an
immediate consequence of Hadamard-Cartan.

\begin{proof}
%
%
Since $C$ is $\varepsilon$-locally connected, for every point $x\in
C$, the set of points of~$C$ which 
can be joined to $x$ by a \emph{rectifiable} curve is a path component
of~$C$, hence all of $C$. So every pair of points of $C$ can be joined
by a rectifiable curve.

For $x,y\in C$, let $l(x,y)$ be the infimum of possible lengths of
curves in $C$ joining $x$ and $y$. 

We argue by induction on $n$ and show: if
$l(x,y)<n\varepsilon$, then $\overline{xy}\subset C$ (and
$l(x,y)=d(x,y)$). For $n=1$, the claim is trivial. 

Assume the claim to be true for $n$, and let $x,y$ be such that
$l(x,y)\in [n\varepsilon ,(n+1)\varepsilon )$. Let
$g_{m}:[0,l(x,y)]\rightarrow C$ be curves of constant speed, such that
$l(g_{m})<(n+1)\varepsilon$ and $l(g_{m})\searrow l(x,y)$.

Let $p_{m}:=g_{m}(t)$ be such that $d(p_{m},y)=\varepsilon$ (such a
point exists; otherwise, the claim were trivial). 
We have $d(x,p_{m})\leq l(x,p_{m})\leq l(g_{m}|_{[0,l(x,y)-\varepsilon]})
<n\varepsilon$, so
by induction hypothesis, we
have $\overline{xp_{m}}\cup \overline{p_{m}y}\subset C$, and we may
assume that $g_{m}$ is a parametrization of these two
segments. Let
$q_{m}:=\overline{p_{m}x}(\varepsilon )$ (as above, $q_{m}$ has to
exist in order for the claim to be non-trivial: if $q_{m}$ does not
exist, then $\{x,y \}\subset B_{\varepsilon}(p_{m})$, so
$\overline{xy}\subset C$ by $\varepsilon$-local-convexity).

We examine the comparison angle $\tilde{\angle}_{p_{m}}(q_{m},y)$:
Since $C$ is $\varepsilon$-locally convex around $p_{m}$, we have
$\overline{q_{m}y}\subset C$. Therefore, the comparison angle has to
be large when $m$ is large:  $d(x,q_{m})+d(q_{m},y)\geq
l(x,y)\swarrow l(g_{m})$, implies $d(q_{m},p_{m})+d(p_{m},y)-
d(q_{m},y) = 2\varepsilon -d(q_{m},y) \rightarrow 0$.

Hence,
$\tilde{\angle}_{p_{m}}(q_{m},y)\rightarrow \pi$. 
Since $q_{m}\in \overline{xp_{m}}$, we have
$\tilde{\angle}_{p_{m}}(q_{m},y)\leq 
\tilde{\angle}_{p_{m}}(x,y)$. So for large $m$, the union
$\overline{xp_{m}}\cup 
\overline{p_{m}y}\subset C$ is almost a geodesic; in particular,
we have $l(x,y)=d(x,y)$ and it is now immediate that the $g_{m}$
converge to $\overline{xy}$, finishing the proof.
%
%
%
%
\end{proof}

\begin{remark}\label{boundaryRemark}
Let $C$ be a closed connected subset of $X$, and
$\partial C$ be the (usual) boundary of $C$ as a topological subset of
$X$. 
Assume that for some $\varepsilon >0$ and every $x\in \partial C$ we
have convexity of 
$B_{\varepsilon}(x)\cap C$. Then $C$ is $\varepsilon /2$-locally
convex, hence convex. 
\end{remark}

Similarly, we have the following lemma:
\begin{lemma}\label{boundaryLemma}
Let $C_{1},C_{2}$ be two closed convex subsets of $X$ and $\varepsilon
>0$. If 
$C_{1}\cup C_{2}$ is connected, and $B_{\varepsilon}(x)\cap (C_{1}\cup
C_{2})$ is convex for every $x\in \partial C_{1}\cap \partial C_{2}$,
then $C_{1}\cup C_{2}$ is convex.
\end{lemma}
\begin{proof}
First, we will show that for every $x\in C_{1}\cap C_{2}$, we have
convexity of $B_{\varepsilon /2}(x)\cap (C_{1}\cup C_{2})$. Then we
show that this is sufficient.

So let $x\in C_{1}\cap C_{2}$. We show directly that for $y,y'\in
B_{\varepsilon /2}(x)\cap (C_{1}\cup C_{2})$, we have
$\overline{yy'}\subset C_{1}\cup C_{2}$. Assume that this is not the
case (for some $x,y,y'$ as above). Then by assumption, we have
$\partial C_{1}\cap \partial C_{2}\cap 
B_{\varepsilon /2}(x)=\emptyset$. 

\Wlo, we have $y\in C_{1}\backslash C_{2}$ and $y'\in C_{2}\backslash
C_{1}$. Let $y_{t}:=\overline{yy'}(t)$, and let $z_{t}$ be the
endpoint of the segment $\overline{xy_{t}}\cap C_{1}$. 
Let $T\geq 0$ be the minimal real number such that
for the interval $(T,d(y,y')]$, we have $z_{t}\not =y_{t}$. For (the
closure of) this
interval, we have $z_{t}\in \partial C_{1}$. 

Similarly,
define $z_{t}'$, and obtain an interval $[0,T')$ such that in this
interval, $z_{t}'\in \partial C_{2}$. By assumption, there is a point
$y_{T''}\not \in C_{1}\cup C_{2}$, so the two intervals introduced
above intersect. 

Consider the functions
$d(x,z_{t})/d(x,y_{t})$ and $d(x,z_{t}')/d(x,y_{t})$ on $[T,T']$. Both are
continuous (observe that $x\not \in \overline{yy'}$),
and by the intermediate value theorem,
they are equal at some point. But then, we 
have found 
a 
point of $\partial C_{1}\cap \partial C_{2}\cap B_{\varepsilon
/2}(x)$, in contradiction to the assumption that $B_{\varepsilon
/2}(x)\cap (C_{1}\cup C_{2})$ is not convex.


Now, we want to show that for any $x\in C_{1}\cup C_{2}$, the set
$B_{\varepsilon /8}(x)\cap (C_{1}\cup C_{2})$ is convex.

Assume that this is not the case for some $x,y,y'$ as above. Note that
by the discussion above, we have $d(x,C_{1}\cap
C_{2})>3\varepsilon/8$. Hence, we have 
\[
\min (d(y,C_{1}\cap C_{2}), d(y',C_{1}\cap C_{2})) >\varepsilon /4 ,
~ \text{and trivially}~ d(y,y')\leq \varepsilon /4.
\]
Let $p:=\pi_{C_{1}\cap C_{2}}(y)$. 
The comparison angle satisfies 
$\tilde{\angle}_{p}(y,y')< \pi /3$. 

Let
$z:=\overline{py}(\varepsilon/8),z':=\overline{py'}(\varepsilon/8)$. 
By the remark about the comparison angle above, $d(z,z')<\varepsilon
/8$ and 
$d(y,q)<d(y,p)$ for every $q\in \overline{zz'}$. 
Note that $z\in C_{1},z'\in C_{2}$.
Convexity of
$B_{\varepsilon /2}(p) \cap (C_{1}\cup C_{2})$ implies that
$\overline{zz'}\cap C_{1}\cap C_{2}\not =\emptyset$. This is a
contradiction to the definition of $p$. 
\end{proof}

We will use the following reformulation quite often:

\begin{corollary}\label{commonCondition}
Let $C_{1},C_{2},K$ be closed convex subsets of $X$, and let $\varepsilon >0$.
Assume that $(C_{1}\cup C_{2})\cap K$ is connected, and that for every
$x\in \partial C_{1}\cap \partial C_{2}\cap K$, we have convexity of
$B_{\varepsilon}(x)\cap K \cap (C_{1}\cup C_{2})$. 

Then $(C_{1}\cup C_{2})\cap K$ is convex.
\end{corollary}

\begin{proof}
Since $K$ is a CAT(0)-space itself, this is just a redraft of the
previous lemma.
\end{proof}

\def\hyph{\penalty0\hskip0pt\relax}
 
\flushbottom

\def\slope{{\rm slope}}
\def\Fold{\mathop{\rm Fold}}
\def\star{\mathop{\rm star}}
\def\supp{\mathop{\rm supp}}
\def\cotan{\mathop {\rm cotan}\nolimits}
\def\Conv{\mathop {\rm Conv}}
\def\M{\mathcal{M}}
\def\F{\mathcal{F}}
\def\C{\mathcal{C}}
\def\EB{Euclidean building}
\def\HS{Hadamard space}

\def\T{\mathcal{T}}
\def\K{\mathcal{K}}

\def\comment#1{}

\def\Wlo{Without loss of generality}
\def\wlo{without loss of generality}

\section[Remarks on Busemann functions]
{Remarks on Busemann functions and horo\-balls in \EB s}
\label{busemannSection}

In this section, we examine general conditions, under which the union
of parts of horoballs is convex.

\paragraph{Setting:}
Let $X$ be a \EB\ without flat de Rham factor, $\eta_{1},\eta_{2}\in
\partial_{T}X$ be two 
boundary points of the same 
type (not necessarily regular), and $p\in X$. We normalize the
corresponding Busemann functions 
$b_{1},b_{2}$ such that $b_{1}(p)=b_{2}(p)=0$. 

Consider $\eta_{1}$ as a
point in
the (spherical) model apartment~$S$.
Let $\alpha >0$ be the maximal angle such that $\angle (\eta_{1},\eta)\leq
\alpha$ implies that $\eta_{1}$ and $\eta$ lie in a common Weyl
chamber of the Coxeter complex~$(S,W)$.
Since $X$ has no flat de Rham factor, we have $\alpha \leq \pi /2$.

For the first two lemmas,
assume there exists a ray
$\overline{p\xi}$ such that 
\[
\angle_{p}(\eta_{1},\xi )=\angle_{p}(\eta_{2},\xi )=\pi.
\]

Since the set of singular points of the Coxeter complex $(S,W)$ is
invariant under the map sending every point to its antipode,
we have:
Whenever $\angle (\xi,\xi')\leq \alpha $ for some $\xi '\in
\partial_{T}X$, then the points $\xi$ and 
$\xi'$ lie in a common Weyl 
chamber of~$\partial_{T}X$.

Note that this implies in particular: If $\xi '\not =\xi$ has the same
type as $\xi$, then $\angle (\xi ,\xi ')\geq 2\cdot \alpha$.

\begin{lemma}\label{equalBusemannFunctions}
Let $\overline{p\xi'}$ be a ray with $\angle_{p}(\xi ,\xi ')\leq
\alpha$. Then 
\[
b_{1}|_{\overline{p\xi '}}=b_{2}|_{\overline{p\xi '}}.
\]
\end{lemma}
\begin{proof}
Note that any Busemann function $b_{\eta}$ is piecewise linear and
convex along any ray
$\overline{p\xi '}$, the 
slope in $x\in \overline{p\xi '}$ being $-\cos(\angle_{x}(\eta ,\xi'))$
(this is well known; it follows from \cite[4.1.2]{kleinerLeeb}).

Now the possible values of $\angle_{x}(\eta_{i},\xi ')$ form a finite
set (determined by the types of $\eta_{i},\xi '$), 
and if $\overrightarrow{x\xi''}$ is of the same type as
$\overrightarrow{x\xi'}$, then  $\angle_{x}(\hat{\xi} 
,\xi '')\geq \angle_{x}(\xi ,\xi ')$ for every antipode $\hat{\xi}$ of
$\eta_{1}$
(if $\overrightarrow{x\xi ''}$ does not lie in a
common Weyl chamber with $\overrightarrow{x\hat{\xi}}$, then
$\angle_{x}(\hat{\xi} ,\xi '')\geq 
\alpha$). 

So $\angle_{p}(\eta_{1},\xi ') = 
\angle_{p}(\eta_{2},\xi ') = \pi -\angle_{p}(\xi ,\xi ')\in [\pi
-\alpha ,\pi ]$
 is maximal. Since 
the slope of $b_{i}$ is increasing 
along $\overline{p\xi '}$, it is constant, and the claim follows from
our assumption $b_{1}(p)=b_{2}(p)$.
\end{proof}

We continue working in the setting introduced above.

\begin{lemma}\label{tubeIntersectsHoroball}
Let $R>0$, and $D\geq \max(R, R/\tan\alpha )$. 
Then
\[
B_{R}(\overline{p\xi})\cap \{b_{1}\leq D \} = 
B_{R}(\overline{p\xi})\cap \{b_{2}\leq D \}.
\]
\end{lemma}
\begin{proof}
Let $x\in B_{R}(\overline{p\xi})$. If $\angle_{p}(x,\xi)\leq \alpha$,
then $b_{1}(x)=b_{2}(x)$ by the previous lemma, so $x$ is either
contained in both sets, or in none of them.

So it suffices to show: if $\angle_{p}(x,\xi)> \alpha$, 
then $b_{i}(x)\leq D$ for both $i$.\newline
Let
$x':=\pi_{\overline{p\xi}}(x)$. 
We may assume $x'\not =p$ because of $D\geq R$.
Consider a point $y\in \overline{xx'}$
such that $\angle_{p}(y,\xi )=\alpha$. Then
$b_{1}(y)=b_{2}(y)=d(p,y)\cdot \cos\alpha$, and
$d(y,\overline{p\xi}) = d(y,x')\geq d(p,y)\cdot \sin\alpha$. We have 
\begin{align*}
b_{i}(x)&\leq b_{i}(y) + (R - d(y,\overline{p\xi}))
\leq R\cdot
\Kr{1+\frac{d(p,y)}{R}\cdot (\cos\alpha -\sin\alpha )}\\
&=R\cdot
\Kr{1+\frac{d(p,y)\cdot \sin\alpha}{R}\cdot (\cotan\alpha -1 )}
\end{align*}
If $\alpha \geq \pi /4$, we have $\cotan \alpha \leq 1$
so the inequality above 
implies $b_{i}(x)\leq R\leq D$.

If $\alpha \leq \pi /4$, we have $\cotan \alpha \geq 1$, and
we use $d(p,y)\leq R/\sin\alpha$: Then the 
inequality above becomes $b_{i}(x)\leq R\cdot (1+\cotan \alpha -1)\leq D$.
\end{proof}

From now on, we do not require the existence
of a common antipode  $\overrightarrow{p\xi }$ of the two
$\overrightarrow{p\eta_{i}}$ in $\Sigma_{p}(X)$ anymore.




%

\begin{lemma}\label{joiningHoroballs}
Let $D>R\cdot \cos\alpha   > 0$. Then the set 
\[
C:=B_{R}(p)\cap \Kr{\{b_{1}\leq D\} \cup \{b_{2}\leq  D\}}
\]
is convex.
\end{lemma}
\begin{proof}
We want to apply Corollary \ref{commonCondition}:
It suffices to find an $\varepsilon >0$ such that 
 for any point $x\in C$ with
$b_{1}(x)=b_{2}(x)=D$, we have convexity of
$B_{\varepsilon}(x)\cap C$.

Let us first choose the $\varepsilon$, depending only on the type of
$D,R$, and $\alpha$ (but not on a specific point $x\in C$):

Let $\delta :=(D-R\cdot \cos\alpha )/2$, and choose $\hat{\alpha}<\alpha$ such
that $R\cdot \cos\hat{\alpha}\leq D-\delta$. 
\newline
Now consider a Euclidean triangle $A,B,C$ with $d(A,B)\geq \delta$ and
$\angle_{A}(B,C)\geq \alpha -\hat{\alpha}$. Let $\varepsilon '$ be such
that $d(B,C)\geq \max(\varepsilon ',\varepsilon '/\tan\alpha )$. Set
$\varepsilon :=\min 
(\delta ,\varepsilon ')$. 

%
%

Now let $x\in C$ be a point with $b_{1}(x)=b_{2}(x)=D$.\newline
There is a finite subdivision $(x_{0}=x, x_{1},\dotsc ,x_{m}=p)$ of
$\overline{xp}$ such that
$\Conv(x_{j},x_{j+1},\eta_{i})$\footnote{Throughout this paper, 
$\Conv$ always denotes the convex hull of its arguments. We use it
with a variety of different kinds of arguments, but no confusion
should arise.} 
is isometric
to a flat half-strip (for $0\leq j<m, 1\leq i\leq 2$) (see
\cite[4.1.2]{kleinerLeeb}). 

Now $D/R>\cos\alpha$ implies that both $b_{i}$ have
maximal slope on the segment~$\overline{x_{1}x}$.

In fact, since $R\cdot \cos\alpha < D- \delta$, we have
$d(x_{1},x)>\delta$ (recall that if the slope of $b_{i}$ is not
maximal, then it is at most $\cos\alpha$).

For a point $y\in \overline{x_{1}x}\backslash \{x_{1} \}$, we have
$\angle_{y}(x_{1},\eta_{i})\leq \hat{\alpha}$ (since the slope of $b_{i}$
along $\overline{x_{1}x}$ has to be larger than $\cos\hat{\alpha}$),
hence $\angle_{y}(\eta_{1},\eta_{2})<2\alpha$, and
$\overrightarrow{y\eta_{1}}=\overrightarrow{y\eta_{2}}$. 

Let $x'\in \overline{x\eta_{1}}\cap \overline{x\eta_{2}}$
be such that $\overrightarrow{x'\eta_{1}}\not
=\overrightarrow{x'\eta_{2}}$. We have
$\angle_{x_{1}}(x',\eta_{i})\leq \pi -\alpha$
(otherwise, we would obtain
$\overrightarrow{x'\eta_{1}}=\overrightarrow{x'\eta_{2}}$ as above), and
$\angle_{x_{1}}(x,\eta_{i})\geq \pi -\hat{\alpha}$. 
This implies that we have $d(x',x)\geq \varepsilon$ by construction.

Of course, $b_{1}(x')=b_{2}(x')=b_{i}(x)-d(x,x') = D - d(x,x')$. By Lemma
\ref{tubeIntersectsHoroball}, we have  
\begin{align*}
B_{\varepsilon}(x)\cap \{b_{1}\leq D \}&= 
B_{\varepsilon }(x)\cap \{b_{2}\leq D\}.
\end{align*}
Hence, $B_{\varepsilon}(x)\cap C$ is convex, and Corollary
\ref{commonCondition} applies. 
\end{proof}

\section{Geometry of \EB s of type $A_{2}$}

In the remainder of this paper, we will work with \EB s of type
$A_{2}$.
Their Tits
boundaries are spherical buildings of type~$A_{2}$.

The (spherical) Coxeter complex $A_{2}$ is the unit circle with the
group of symmetries of an isosceles triangle acting on it (see Figure
\ref{boundaryFigure}). A (discrete) Euclidean Coxeter complex of type
$A_{2}$ is the Euclidean plane, tesselated by isosceles triangles (see
Figure \ref{A2-figure}).

The most important example of a \EB~of type $A_{2}$ is the building
associated to $SL(3,\Q _p)$; its geometry is described in detail in
\cite{kremserDipl}. 

\subsection{The spherical building structure of $\partial_{T}X$}
Let $X$ be a \EB~of type $A_{2}$. Then the boundary at infinity
$\partial_{T}X$ carries the structure of a spherical
building of type $A_{2}$. Similarly, the space of directions
$\Sigma_{x}(X)$ for any $x\in X$ carries such a structure as well.

Every apartment in such a spherical building $B$ consists of six Weyl chambers of
length $\pi /3$. The vertices 
(the singular points of $B$, the ends of the Weyl
chambers) have two different types (see Figure \ref{boundaryFigure}).

\subsection{Holonomy in spaces modeled on $A_{2}$}
Let $\eta\in \partial_{T}X$ be a regular boundary point. By the
remarks above, $X_{\eta}\simeq \R$. Since every apartment asymptotic
to $\eta$ represents all strong asymptote classes, we get an
orientation on $X_{\eta}$ (induced from a choice of orientation on 
Weyl chambers, determined by the two types of boundary points). 
Then every holonomy map $h_{\eta ,\xi}$ is orientation preserving, and
so is the composition 
\[
h_{\eta_{1},\eta_{2},\eta_{3}}:=h_{\eta_{3},\eta_{1}}\circ
h_{\eta_{2},\eta_{3}}\circ h_{\eta_{1},\eta_{2}}:
X_{\eta_{1}}\rightarrow X_{\eta_{1}} 
\]
for any triple $\eta_{1},\eta_{2},\eta_{3}\in \partial_{T}X$ of
pairwise antipodal regular boundary points.
Such a holonomy map, as an orientation preserving isometry of $\R$, is
just a translation. We will call the translation length of such a
triple its \emph{shift}.

Observe that for the other three 2-dimensional Coxeter complexes,
holonomy maps are orientation-reversing. This major difference
makes it hard to predict a general (2-dimensional) version of
 Theorem \ref{convBuild}.

\subsection{Geometric lemmas for buildings of type $A_{2}$}
\label{geometrySection}

In this section, we collect some geometric properties of \EB s of type
$A_{2}$ that will be useful later. 

\begin{lemma}\label{allAnglesMaximal}
Let $X$ be a building of type $A_{2}$, $p\in X$ and
$\eta_{1},\eta_{2},\eta_{3}$ be 
three singular boundary points of the same type, such that the
$\overrightarrow{p\eta_{i}}$ 
span a flat in $\Sigma_{p}(X)$. Normalize such that $b_{i}(p)=0$. Let
$q\in X$ be 
such that $\angle_{p}(q,\eta_{i})=2\pi /3$ for all $i$ and
$R:=d(q,p)>0$ (so $b_{i}(q)=R/2$). 

Let 
$C:=\{x\st \text{at least two }b_{i}(x)\leq R/2 \}$.
Then $K:=C\cap 
B_{R\cdot \sqrt{3}/2}(q)$ is convex.
More specifically, there exist $i,j$ such that $K=B_{R\cdot
\sqrt{3}/2}(q)\cap \{b_{i}\leq R/2 \}\cap \{b_{j}\leq R/2 \}$. 
\end{lemma}
\begin{proof}
%
%
%
Pick an $x\in B_{R\cdot \sqrt{3}/2}(q)$, and observe that
$\angle_{p}(x,q)\leq \tilde{\angle}_{p}(x,q)\leq \pi /3$ (by triangle
comparison).  

We distinguish two cases: The first case is that the initial
directions of
$\overline{\overrightarrow{pq}\overrightarrow{p\eta_{i}}} \subset
\Sigma_{p}X$ are all 
distinct.
Then, there are two $i$ such that
$\angle_{p}(x,\eta_{i})\geq 2\pi /3$. So $x\in C$ if and only if all
three Busemann functions are at most $R/2$. (In this case, we can
choose $i,j$ arbitrarily.)

Otherwise, (exactly) two of the above-mentioned initial directions
agree (\wlo, those corresponding to $1,2$; these correspond to the
$i,j$ in the claim). 

We claim that $K=B_{R\cdot \sqrt{3}/2}(q)\cap \{b_{1}\leq R/2 \}\cap
\{b_{2}\leq R/2 \}$.\newline
 Indeed, if $b_{3}(x)\leq R/2$ for $x\in K$, then either
$\angle_{p}(x,\eta_{3})\geq 2\pi /3$ and
$b_{3}(x)\leq \min (b_{1}(x),b_{2}(x))$, or
$\angle_{p}(x,\eta_{1})=\angle_{p}(x,\eta_{2})> 2\pi /3$. In the last
case, we have $b_{1}(x)=b_{2}(x)$, hence the claim follows.
%
%
\end{proof}

\bigbreak

\begin{figure}\def\color[#1]#2{}\def\mddefault{}
\def\updefault{}
\vspace*{2.7cm}
\hfill 
\input{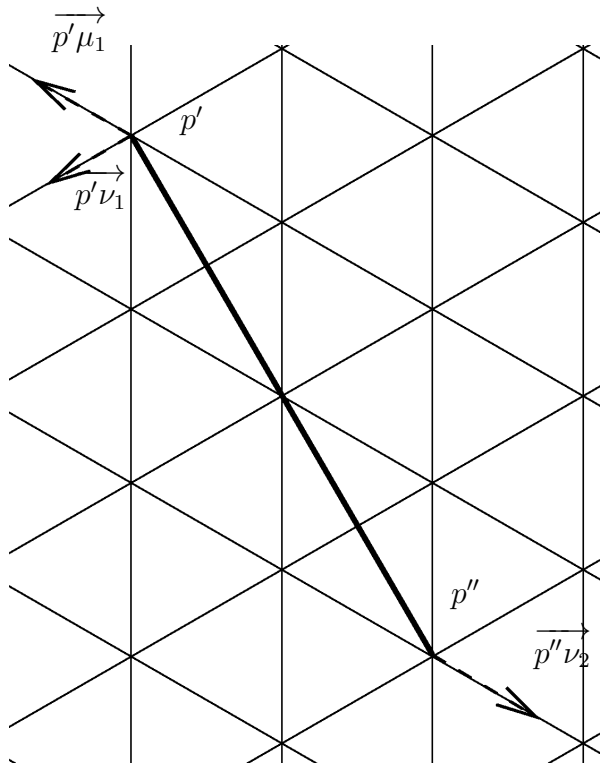}
  \hfill ~
\caption{The setting of the lemmas}
\label{A2-figure}
\end{figure}

For the next two lemmas, we need a setting which will be introduced in
more detail later:

Let $\eta_{1},\eta_{2}$ be antipodal centers of Weyl chambers in the
boundary of a \EB \ $X$ of type $A_{2}$, and let $F_{1,2}$ be the flat
containing $\eta_{1},\eta_{2}$ in its boundary. Let
$\nu_{1},\eta_{1,2},\nu_{2},\mu_{2},\xi_{1,2},\mu_{1}$ be the
singular points in  
$\partial_{T} F_{1,2}$ as in Figure \ref{boundaryFigure}. Figure
\ref{A2-figure} shows a part of the flat $F_{1,2}$, with the boundary
being aligned as in Figure \ref{boundaryFigure} (with $i=1$, $j=2$).

\begin{lemma}\label{triangleContradiction}
Let $p',p''\in F_{1,2}$ such that $b_{1}(p')\leq b_{1}(p'')$.

Let $0\leq \hat{\alpha}<\pi /3$. Then
there is no  $x\in X$ 
with 
\begin{align*}
\angle_{p'}(x,\nu_{1})\,&< \pi
/3+\hat{\alpha},
 &\angle_{p''}(x,\nu_{2})\,&<\pi /3+\hat{\alpha},\\
\text{but }
\angle_{p'}(x,\eta_{1,2})\,&>\hat{\alpha},&
\angle_{p''}(x,\eta_{1,2})\,&>\hat{\alpha}.
\end{align*}
\end{lemma}
\begin{proof}
If $p'=p''$, there is nothing to show.

\Wlo, we assume $b_{1,2}(p'')\leq b_{1,2}(p')$. \newline
Let $\beta :=
\angle_{p''}(p',\xi_{1,2}) = \angle_{p'}(p'',\eta_{1,2})\leq \pi
/2$. 

If $\pi /3\leq \beta\,\, (\leq \pi /2)$, one 
obtains a contradiction to the sum of angles in a triangle:
Indeed, we have 
\begin{align*}
\angle_{p''}(p',x)> &\,\min (\pi +\hat{\alpha}-\beta , \pi /3+\beta
-\hat{\alpha})\\ 
\angle_{p'}(p'',x)> &\,\min (\beta +\hat{\alpha},
\frac{4\pi}{3}-\beta -\hat{\alpha}).
\end{align*}
We see that if $\beta \geq \pi /3$, the sum of these two angles is
greater than $\pi$. Therefore, we have  $\beta <\pi /3$. 
Then
\begin{align}\label{fangle}
\angle_{p''}(p',x)>\pi /3 + \beta -\hat{\alpha},
\end{align}
and 
\begin{align}\label{angleIneq}
\angle_{p'}(p'',x)>\hat{\alpha}+\angle_{p'}(p'',\eta_{1,2}) = 
\hat{\alpha}+\beta  
\end{align}

Let $(p_{0}=p'',p_{1},\dotsc, p_{n}=p')$ be a finite subdivision of
$\overline{p''p'}$ such that each triangle $\Delta (p_{i},p_{i+1},x)$
is flat.

Let $i_{0}>0$ be such that 
the
initial directions of
$\overline{\overrightarrow{p_{i_0}p''}\overrightarrow{p_{i_0}x}}$ and
$\overline{\overrightarrow{p_{i_0}p''}\overrightarrow{p_{i_0}\nu_{2}}}$
 agree $(*)$. 

We will show by induction that every $1\leq i_{0}\leq n$ has this
property, and obtain a contradiction for $i_{0}=n$.

\paragraph{Base case:}
$i_{0}=1$ has Property $(*)$.

If this is not the case, then 
the starting direction of
$\overline{\overrightarrow{p''p_{1}}\overrightarrow{p''x}}$ has to be
different from the starting direction of
$\overline{\overrightarrow{p''p_{1}}\overrightarrow{p''\nu_{2}}}$
(since the triangle $\Delta (p'',p_{1},x)$ is flat). If
this is the case, we have
$\angle_{p''}(p',x) = \angle_{p''}(p_{1},x)>\pi
-(\hat{\alpha}+\beta )$. This is a contradiction to \eqref{angleIneq}.

\paragraph{Claim:}
If $i_{0}<n$ has property $(*)$, then
the initial directions of
$\overline{\overrightarrow{p_{i_0}p'}\overrightarrow{p_{i_0}\nu_{2}}}$
and $\overline{\overrightarrow{p_{i_0}p'}\overrightarrow{p_{i_0}x}}$
agree as well.

Observe that
$\angle_{p_{i_{0}}}(p'',x)<2\pi /3+\hat{\alpha}-\beta$ (by
\eqref{fangle}). 
Assume that the claim is false: 
Then 
\[
\angle_{p_{i_{0}}}(p',x)>\pi -
(\underbrace{(\pi /3+\hat{\alpha})}_{>
\angle_{p_{i_{0}}}(\nu_{2},x)}-\underbrace{(\pi /3-\beta 
)}_{\angle_{p_{i_{0}}}(p'',\nu_{2})}) = \pi 
-\hat{\alpha}+\beta.
\]
Together with \eqref{angleIneq}, this is a
contradiction. 

Now this claim implies $(*)$ for $i_{0}+1$ (by the same argument as in
the base case, with $p_{i_{0}}$ taking the place of $p''$), and it
follows by induction that property~$(*)$ holds for all $1\leq 
i_{0}\leq n$. 

For $i_{0}=n$, we get $\angle_{p'}(p'',x)> \pi
-\hat{\alpha}+(\pi /3-\beta ) = 4\pi/3 -(\hat{\alpha}+\beta)$. 

If this is less than $\pi$, we continue our calculation:
\[
\angle_{p'}(p'',x)+\angle_{p''}(p',x)>5\pi /3-2\hat{\alpha}.
\]
This is a contradiction, since $\hat{\alpha}<\pi /3$.
\end{proof}

\begin{remark}\label{triangleContradictionRemark}
In the statement of the lemma, we can replace the directions
$\overrightarrow{p'\nu_{1}}$ and 
$\overrightarrow{p''\nu_{2}}$ by any other directions antipodal to
$\overrightarrow{p'\mu_{2}},
\overrightarrow{p''\mu_{1}}$ resp. Of course, we also have to adjust
the assumptions after the ``but''. We will usually take care of  this
by showing that $\angle_{p''}(\xi_{1,2},x)<\pi -\hat{\alpha}$ and 
$\angle_{p'}(\xi_{1,2},x)<\pi -\hat{\alpha}$.
\end{remark}

\begin{lemma}\label{triangleContradictionTwo}
Let $p',p''\in F_{1,2}$ such that $b_{1}(p')\leq b_{1}(p'')$.

Assume that
$\angle_{p''}(p',\eta_{1,2})\geq \pi /3$. Let $0\leq \hat{\alpha}\leq \pi
/6$. Then there is no $x\in X$ 
with 
\begin{align*}
\angle_{p'}(x,\mu_{1})\,&< \pi
/3+\hat{\alpha},&\angle_{p''}(x,\nu_{2})\,&<\pi /3+\hat{\alpha}\\
\text{but }
\angle_{p'}(x,\xi_{1,2})\,&>\hat{\alpha},&
\angle_{p''}(x,\eta_{1,2})\,&>\hat{\alpha}.
\end{align*}
\end{lemma}
\begin{proof}
As above, we may assume $p'\not =p''$.

We distinguish two cases: The first case is
$\angle_{p''}(p',\xi_{1,2})=:\beta \in [\pi /3,2\pi /3]$. 
In this case, we have 
\begin{align*}
\angle_{p''}(p',x)\,&>\min (\pi -\beta +\hat{\alpha}, 
\underbrace{\beta +\pi
/3-\hat{\alpha}}_{\geq \pi /2})\\
\angle_{p'}(p'',x)\,&>\min (\pi -\beta +\hat{\alpha}, 
\underbrace{\beta +\pi
/3-\hat{\alpha}}_{\geq \pi /2}).
\end{align*}
Adding these angles, the only case in which we do not get a
contradiction to the sum of angles in a triangle is, if
$\overrightarrow{p''x}\in
\overline{\overrightarrow{p''\eta_{1,2}}\overrightarrow{p''\nu_{2}}}$,
$\overrightarrow{p'x}\in
\overline{\overrightarrow{p'\xi_{1,2}}\overrightarrow{p'\mu_{1}}}$ and
$\beta >\pi /2+\hat{\alpha}$. 

Now this case can be finished as in the lemma above: 
The deciding
inequalities are 
\begin{align*}
\angle_{p_{i_{0}}}(p'',x)&<\beta -\hat{\alpha},\qquad
\angle_{p_{i_{0}}}(p',x)<\beta -\hat{\alpha},
\\
\angle_{p_{i_{0}}}(p'',\nu_{2})&=\beta -\pi /3.
\end{align*}

The second case is $\beta <\pi /3$.
Now we have $\angle_{p'}(p'',x)>\beta +(\pi /3-\hat{\alpha})$ and
$\angle_{p''}(p',x) >(\pi /3-\hat{\alpha})$.
Again, we have to have $\overrightarrow{p''x}\in
\overline{\overrightarrow{p''p'}\overrightarrow{p''\nu_{2}}}$ and
$\overrightarrow{p'x}\in
\overline{\overrightarrow{p'p''}\overrightarrow{p'\mu_{1}}}$.  

To be able to ``switch sides'', we would need a $p_{i_{0}}$ with
$\angle_{p_{i_{0}}}(p'',x)$ at least $\underbrace{\pi
/3-\beta}_{\angle_{p_{i_{0}}}(p'',\nu_{2})} +2\pi /3 
-\hat{\alpha} = \pi -\beta -\hat{\alpha}$, which is impossible
(because $\hat{\alpha}\leq \pi /6$).  
\end{proof}

\begin{remark}
Again, this lemma remains true if we replace
$\overrightarrow{p''\nu_{2}}$ and/or $\overrightarrow{p'\mu_{1}}$ 
by other directions antipodal to
$\overrightarrow{p''\mu_{1}}, \overrightarrow{p'\nu_{2}}$ resp. (and
again, we also have to adjust the assumptions after the ``but'').
\end{remark}

\section{Necessary conditions: S-sets}
Now we shift gears, and turn directly to the proof of Theorem
\ref{convBuild}. We start by examining necessary conditions and
obtaining more and more structure on the sets satisfying the (obvious)
necessary conditions.

\label{necessarySection}
Let us first state the precise definition of a convex rank 1-set:

\begin{definition}\label{convexRank1Definition}
A subset $C\subset X$ of a \HS \ $X$ is called \emph{convex rank 1},
if it is closed, convex, has at least 3 boundary points at infinity
and satisfies: $\partial_{T}C$ is a 
0-dimensional building (i.e.: for all $ \eta ,\xi\in \partial_{T}C$
with $\eta \not =\xi$, we have $\angle_{T}(\eta ,\xi )\geq \pi$).
\end{definition}

Observe that the restriction $|\partial_{T}C|\geq 3$ is not serious:
Every pair of antipodal points in $\partial_{T}X$ (for $X$ a symmetric
space 
or a 
\EB ) can be joined by a geodesic.

From now on, we  focus on a special class of buildings:
In the remainder of this article, $X$ will always stand for a building
of type $A_{2}$. 

In this section we examine necessary conditions for points
$\eta_{i}\in \partial_{T}X$ to be in the boundary of a convex rank 1
set. The most important necessary condition is that there has to be a
tripod for every triple of asymptotic boundary points (Proposition
\ref{existenceOfTripods}). 
We also examine the structure of the set of singular points of
these tripods, and we will obtain a metric tree which is
closely related.

\begin{lemma}\label{boundaryPointsCenters}
If there are at least  three points $\eta_{i}$, then to be pairwise
antipodal, it is necessary that each $\eta_{i}$ is the center of a
Weyl chamber.
\end{lemma}
\begin{proof}
In the Coxeter complex $A_{2}$, the centers of Weyl chambers are the
only points which have the following property: An antipode has the same
type. This property is necessary, since the points $\eta_{i}$ have to
be pairwise antipodal.
\end{proof}

There is another 
obvious necessary condition: Let $A=\partial_{T}C$ be the asymptotic
boundary of a convex rank 1-set. 
Consider a triple of boundary points. Then the corresponding holonomy
map has to have a fixed point (see Section~\ref{holonomySection});
otherwise, every convex set containing the given triple in its
boundary contains a half-plane, and hence does not have rank 1.

Since our boundary points are regular,
the holonomy map of a triple is an isometry of $\R$ to itself. This map is also
orientation preserving, so it is just a translation, 
determined by its translation length, which we call its shift.

So the necessary condition is:
For each triple of points of
$A$, their shift has to be 0 (i.e.~the holonomy map has to be the
identity map).

\subsection{Notation}\label{notationSection}

\begin{figure}\def\color[#1]#2{}\def\mddefault{}
\def\updefault{}

\hfill 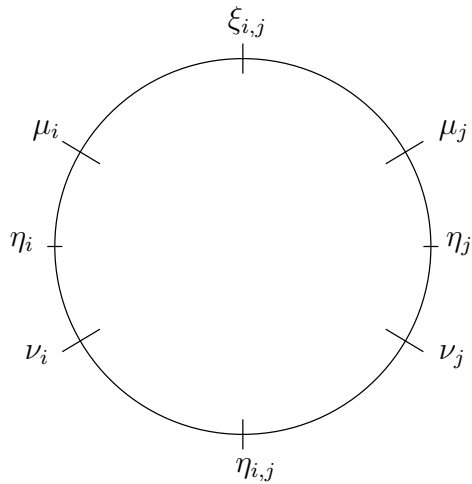\hfill ~
\caption{The apartment $\partial_{T}F_{i,j}\subset \partial_{T}X$ with
its singular points
$\nu_{i},\eta_{i,j},\nu_{j},\mu_{j},\xi_{i,j},\mu_{i}$ and the two
regular points $\eta_{i},\eta_{j}$.}
\label{boundaryFigure}
\end{figure}

\begin{definition}
A subset $A\subset \partial_{T}X$ with $|A|\geq 3$ is an \emph{S-set},
if the points of 
$A$ are pairwise antipodal (i.e.\ $A$ is a 0-dimensional subbuilding),
and for each 
triple of points of $A$, the shift is 0.
\end{definition}

In what follows, $A=\{\eta_{i}\st i\in I \}$ will always be an S-set
(see also the definition of a \emph{good S-set}, \ref{goodDef}).  

A \emph{tripod} is a metric tree with three asymptotic boundary
points. It may  also be viewed as the Euclidean cone over a set of
cardinality three. A tripod in $X$ is determined by a (singular) point
$p$ and three boundary points $\xi ,\nu ,\mu $. This data determines a
tripod $(p,\xi ,\nu ,\mu )=\Conv (p,\xi ,\nu ,\mu )$ if and only if
$\angle_{p}(\xi ,\nu )=\angle_{p}(\nu ,\mu 
)=\angle_{p}(\mu ,\xi )=\pi$. 

In our setting, a \emph{tripodal point}
$p_{i',j',k'}$ (for 
three distinct $i',j',k'\in I$) is a point such 
that $(p_{i',j',k'},\eta_{i'},\eta_{j'},\eta_{k'})$ determines a tripod. 
When a tripodal point is given, then $T_{i,j,k}$ denotes the corresponding
tripod. If the tripodal point is to be emphasized, we also say that 
$(p_{i',j',k'},\eta_{i'},\eta_{j'},\eta_{k'})\in
X\times (\partial_{T}X)^{3}$ is a tripod.

For $i\in I$, let $\nu_{i},\mu_{i}$ be the endpoints of the Weyl
chamber spanned by 
$\eta_{i}$, such that all the $\nu_{i}$ have the same type.

For a pair $i,j\in I$, let $\eta_{i,j}$ be the center of the geodesic
$\overline{\nu_{i}\nu_{j}}\subset \partial_{T}X$; similarly define
$\xi_{i,j}$ (see Figure \ref{boundaryFigure}). 
Let $F_{i,j}$ denote the unique flat in $X$
such that $\eta_{i},\eta_{j}\in \partial_{T}F_{i,j}$; so the singular
vertices of $\partial_{T}F_{i,j}$ are
$\nu_{i},\eta_{i,j},\nu_{j},\mu_{j},\xi_{i,j},\mu_{i}$. 

For a triple $i,j,k\in I$, let $l_{i,j,k}:=F_{i,j}\cap F_{j,k}\cap
F_{i,k}$. 
By definition of tripods, $l_{i,j,k}$ is precisely the set of tripodal
points for 
$(\eta_{i},\eta_{j},\eta_{k})$. 
It follows from the flat strip theorem
(\cite[II.2.13]{bridsonHaefliger}), that for every pair of 
tripods for a given triple of boundary points, they are parallel
to each other, and their convex hull splits as a product ``tripod
$\times$ interval''. 

We will see below that $l_{i,j,k}$ is a
non-empty line
segment (which may degenerate to a point, a ray or a geodesic line).
We will say that $p$ is the ``lower endpoint'' of $l_{i,j,k}$ if $p$
minimizes $b_{i,j}|_{l_{i,j,k}}$ (it follows from Proposition
\ref{existenceOfTripods} that such a point
exists). Analogously, we define the ``upper endpoint'' of
$l_{i,j,k}$. 

Let us introduce another last piece of notation:
When we  deal with the points
$\eta_{i},\eta_{i,j},\xi_{i,j}\in \partial_{T}X$,
we  simplify notation and call the corresponding Busemann
functions $b_{i},b_{i,j},b_{i,j}'$ respectively (instead of
the standard notation
$b_{\eta_{i}}$, $b_{\eta_{i,j}}$, $b_{\xi_{i,j}}$).

\subsection{Existence of tripods}
\begin{proposition}\label{existenceOfTripods}
Let $X$ be a \EB ~of type $A_{2}$, and let
$\eta_{1},\eta_{2},\eta_{3}\in \partial_{T}X$ be three pairwise
antipodal points. If their shift is 0, then there exists a tripod
$(p,\eta_{1},\eta_{2},\eta_{3})$. 
\end{proposition}

Note that the proposition can also be phrased as follows:
Every S-set of cardinality 3 is the asymptotic boundary of a
convex rank 1-set, and this rank~1-set can be chosen to be a tripod.
\begin{figure}\def\color[#1]#2{}\def\mddefault{}
\def\updefault{}

\hfill \includegraphics{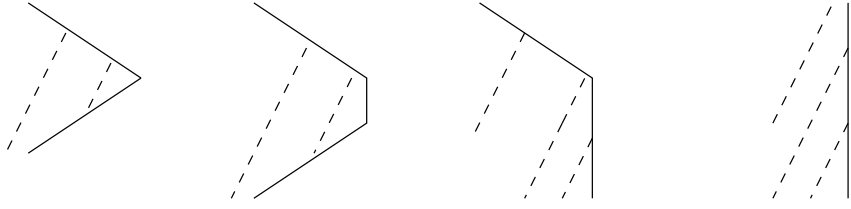}\hfill ~
\caption{the possibilities for $F_{1,2}\cap F_{1,3}$ (where $\eta_{1}$
is lying ``on the left'')}
\label{flatsIntersectFigure}
\end{figure}
\begin{proof}
Observe that $F_{1,2}\cap F_{1,3}=:S$ is non-empty, closed and convex
(by \cite[4.6.3]{kleinerLeeb}).
The Busemann function $b_{1}$ is bounded above on $S$, since otherwise
we have $\angle_{T}(\eta_{2},\eta_{3})<\pi$ (note that $\partial S$
is a polygonal curve consisting of at most three line
segments/rays/lines, 
see Figure \ref{flatsIntersectFigure}).

Let $p$ be an extremal point for $b_{1}|_{S}$. We claim that
$(p,\eta_{1},\eta_{2},\eta_{3})$ is a tripod.

Assume that this is not the case. Then
$\angle_{p}(\eta_{2},\eta_{3})\geq \pi 
/3$ (since both directions are the centers of a Weyl chamber, the smallest
non-zero value for their angle is $\pi /3$). 

Since the shift is zero, we obtain points $p',p''$ in
$\overline{p\eta_{2}},\overline{p\eta_{3}}$ resp.~such that
$\overline{p'\eta_{2}}\cup \overline{p'p''}\cup
\overline{p''\eta_{3}}$ is a geodesic line. 

Let us choose $p',p''$ (as $p$) such that
$\angle_{p'}(\eta_{1},\eta_{3})\not 
=0\not =\angle_{p''}(\eta_{1},\eta_{2})$. 
Then each of these angles is at least $\pi /3$, so $\Delta (p,p',p'')$
is a flat isosceles triangle (by triangle rigidity in CAT(0)-spaces, see
\cite[II.2.9]{bridsonHaefliger}).  

Let $\nu$ be the midpoint of the geodesic
$\overline{\overrightarrow{pp'}\overrightarrow{pp''}}\subset
\Sigma_{p}(X)$. Then $\angle_{p}(\eta_{2},\nu) =
\angle_{p}(\eta_{3},\nu )=\pi /6$; 
Observe that $\mu \in \Sigma_{p}(X),\angle_{p}(\eta_{j},\mu )\leq \pi
/6$ implies $\mu \in \Sigma_{p}(F_{1,j})$ for $j\in \{2,3 \}$ (because
$\mu$ and $\overrightarrow{p\eta_{j}}$ lie in a common Weyl chamber of
$\Sigma_{p}(X)$).
So  
either
\[
\nu =\overrightarrow{p\nu_{2}} = \overrightarrow{p\nu_{3}}
\text{ or }
\nu =\overrightarrow{p\mu_{2}}= \overrightarrow{p\mu_{3}}.
\]
In both cases, we have
$\overrightarrow{p\nu}\in \Sigma_{p}(F_{1,2})\cap
\Sigma_{p}(F_{1,3}) = \Sigma_{p}(F_{1,2}\cap F_{1,3})$. Since
$\angle_{p}(\eta_{1},\nu )=2\pi /3$, this is a 
contradiction to the construction of $p$.
\end{proof}


The proof also shows that a convex rank 1-subset of $X$ has to contain a
tripod for every triple of boundary points (because for a strong
asymptote class which does not correspond to a tripod, we obtain a
flat isosceles triangle in the convex hull, and its center is a
tripodal point). 

Hence, the following condition is also necessary for an S-set $A$ to
be in the asymptotic boundary of a convex rank 1-set:

\begin{definition}\label{goodDef}
An S-set $A$ is called \emph{good}, if it satisfies the following
condition: 
We can choose tripodal points $p_{i,j,k}$ (for every triple $i,j,k\in
I$) such that 
for all $i'\in I$,
the convex hull of (all) the strong asymptote
classes $[\overline{p_{i',j,k}\eta_{i'}}]$ is bounded.
\end{definition}

\begin{example}\label{fourPointExample}
Let us give an example of a 4-point S-set which does not lie in the
boundary of an embedded tree:

We start with two antipodal centers of Weyl chambers,
$\eta_{1},\eta_{2}$, and pick a singular vertex $p$ in $F_{1,2}$.
Choose a ray $\overline{p\eta_{3}}$, such that
$(p,\eta_{1},\eta_{2},\eta_{3})$ is a tripod.
Let us choose $\eta_{3}$ such that the intersection $F_{1,2}\cap F_{2,3}$
is a flat sector (this corresponds to the left-most set drawn in
Figure \ref{flatsIntersectFigure}).%
\footnote{This is possible in ``most'' Euclidean buildings of type
$A_{2}$; pick the building associated to $SL(3,\Q _5)$
for example.} 

Now pick an inner point $p'$ of $F_{1,2}\cap F_{2,3}$ satisfying
$b_{1,2}(p')\not =b_{1,2}(p)$. As above, pick 
a ray $\overline{p'\eta_{4}}$, such that
$(p',\eta_{1},\eta_{2},\eta_{4})$ is a tripod and $F_{1,2}\cap
F_{1,4}$ is a flat sector. 
By construction, we
have 
\[
\overrightarrow{p'\eta_{1}}=\overrightarrow{p'\eta_{3}},
\]
so we
also have a tripod $(p',\eta_{3},\eta_{2},\eta_{4})$. Similarly, our
construction implies that $p$ lies in the interior of $F_{1,2}\cap
F_{1,4}$, so we also have the tripod
$(p,\eta_{1},\eta_{3},\eta_{4})$. Our choices imply that $p,p'$ are
the unique tripodal points (at least when considered as
$p_{1,2,3},p_{1,2,4}$ resp.), so $b_{1,2}(p)\not =b_{1,2}(p')$ implies
that there is no embedded tree with the given four boundary points.

A similar situation is depicted in Figure \ref{fourPointsFigure}; in
the next section, we are going to show that the general situation is
always similar to the one described here. 

Applying the construction above to obtain an S-set with infinitely
many boundary points, we see that an S-set needs not be good.
\end{example}
\subsection[S-sets with 4 points]{S-sets with 4 points: relative
position of their tripodal points} 

In this section, we examine S-sets $A$ of cardinality 4: 
We show that we can always do with at most 2 tripodal points: If there
is no $4$-pod (i.e.~a Euclidean cone over $A$) embedded in $X$, then
we construct two points, each of which is tripodal for two triples of
points of $A$. 

All of the Lemmas in this section are formulated such that the
assumptions rule out existence of a $4$-pod; only Proposition
\ref{fourSingularPoints} is formulated to make 
sense even in this case.

We also discuss the possible choices for the tripodal points in
question, and the relative position of the two (sets of) points to
each other. These are technical results needed in the sequel.

\begin{lemma}\label{structureLemma}
Let $\{\eta_{1},\eta_{2},\eta_{3},\eta_{4} \}\subset \partial_{T}X$ be
an S-set of cardinality 4. Assume there are tripods
$(\bar{p},\eta_{1},\eta_{2},\eta_{3})$ and
$(\bar{p}',\eta_{1},\eta_{2},\eta_{4})$ with
$b_{1}(\bar{p})<b_{1}(\bar{p}')$. Then there are points $p,p'\in X$
such that we have tripods 
\[
({p},\eta_{1},\eta_{2},\eta_{3}),
({p},\eta_{1},\eta_{3},\eta_{4}), \text{ and }
({p}',\eta_{1},\eta_{2},\eta_{4}),
({p}',\eta_{2},\eta_{3},\eta_{4}).
\] 
In particular,
\[
\overline{pp'}\subset F_{1,2}\cap F_{2,3}\cap F_{3,4}\cap F_{1,4}.
\]
and $\angle_{p}(\eta_{1,2},p')\in [\pi /3,2\pi /3]$. 
\end{lemma}

\begin{proof}
Let us choose tripodal points $p\in l_{1,2,3},p'\in l_{1,2,4}$
such that $d(p,p')$ is minimal; note that $b_{1}(p')-b_{1}(p) =
b_{1}(\bar{p}')-b_{1}(\bar{p})>0 $.
Note that this implies in particular that there is no 4-pod with the
given four boundary points embedded in $X$.

\begin{figure}\def\color[#1]#2{}\def\mddefault{}
\def\updefault{}

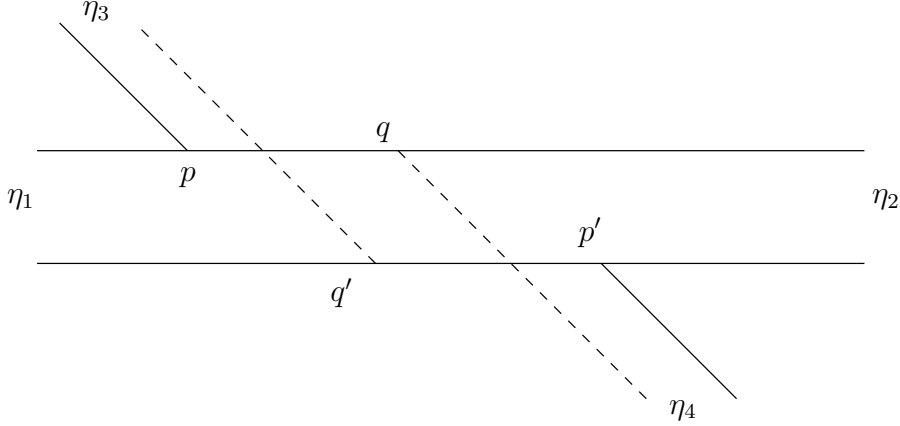
\caption{The situation from Lemma \ref{structureLemma}. Observe that
 $q,q'$ are not tripodal points (unless $q=p'$ and $q'=p$).}
\label{fourPointsFigure}
\end{figure}

If $b_{1,2}(p)=b_{1,2}(p')$, then we have found an isometrically
embedded tree having $\eta_{1},\eta_{2},\eta_{3},\eta_{4}$ as
asymptotic boundary points. The claim of the lemma is now trivial.

So we may assume $b_{1,2}(p)\not =b_{1,2}(p')$, and
\wlo~that
 $b_{1,2}(p)> b_{1,2}(p')$
(by exchanging the $\eta_{i,j}$ and the $\xi_{i,j}$, if necessary);
note that under these assumptions, $p$ is the lower
endpoint of $l_{1,2,3}$, and $p'$ is the upper endpoint of
$l_{1,2,4}$. 
We normalize such that $b_{1,2}(p)=b_{1}(p)=b_{2}(p)=0$.

First, we want to show $p'\in F_{2,3}$. Assume that this is not the
case. 

In $F_{1,2}$, consider the line $l_{1,2}$ 
passing through $p$ with endpoints $\mu_{1},\nu_{2}$. Then the ray
$l_{1,2}\cap \{b_{2}\leq 0 \}$ is a boundary segment of $F_{1,2}\cap
F_{2,3} (\dagger)$. 

Similarly, consider the line $l_{1,2}'$ 
passing through $p'$ with endpoints $\nu_{1},\mu_{2}$.
Then the ray $l'_{1,2}\cap \{b_{2}\leq b_{2}(p') \}$ is a boundary
segment of $F_{1,2}\cap F_{2,4}(\ddagger)$.

The two lines $l_{1,2},l'_{1,2}$ bound a sector $S\subset F_{1,2}$
with tip $p''$, 
containing 
$\eta_{2}$ in its asymptotic boundary. Let $\rho \subset
l_{1,2},\rho'\subset l'_{1,2}$ be its 
bounding rays.

We are assuming that $p'\not \in F_{1,2}\cap F_{2,3}$. 
Since $b_{1,2}(p)> b_{1,2}(p')$, this implies that $p'$ lies
``below" $l$ (otherwise, $p'\in \Conv (p,\eta_{2},\nu_{2})\subset
F_{1,2}\cap F_{2,3}$), see Figure \ref{sector}. 

\textbf{In this case, we claim}
 $S= F_{2,3}\cap F_{2,4}$: The relation $\subset$ is clear (because
$\{p'' \}= l_{1,2}\cap l_{1,2}' \subset  F_{1,2}\cap F_{2,3}\cap F_{2,4}$
by $\dagger$ and $\ddagger$,
and $\overline{\nu_{2}\mu_{2}}\subset \partial_{T} F_{2,j}$ for all
$j$).   
For the other inclusion, observe:
near~$\rho '$, the flat $F_{2,3}$
agrees with $F_{1,2}$, while this is not true for
$F_{2,4}$. Similarly near~$\rho$, the flat $F_{2,4}$ agrees with
$F_{1,2}$, but the flat $F_{2,3}$ does not. 

So $S= F_{2,3}\cap F_{2,4}$ as claimed.
Then $(p'',\eta_{2},\eta_{3},\eta_{4})$
is a tripod by our assumptions and our discussion showing existence of
tripods.

\begin{figure}\def\color[#1]#2{}\def\mddefault{}
\def\updefault{}

\hfill 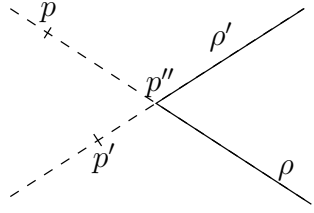\hfill ~
\caption{The relative position of $p,p',p''$.}
\label{sector}
\end{figure}

However, we see immediately that $\angle_{p''}(\eta_{3},\eta_{4})\leq
2\pi /3$: Indeed, we have 
\begin{align*}
\angle_{p''}(\eta_{3},p)\,&=\pi /6 \text{ (by $\dagger$, we have
$\overrightarrow{p''p}\in \Sigma_{p''}( F_{2,3})$)}, \\ 
\angle_{p''}(p,p')\,&=\pi /3, \\
\angle_{p''}(p',\eta_{4})\,&=\pi /6\text{ (by $\ddagger$, we have
$\overrightarrow{p''p'}\in \Sigma_{p''}(F_{2,4})$)}. 
\end{align*}
This contradiction shows $p'\in F_{1,2}\cap F_{2,3}$. 
At the same time, this shows 
$\angle_{p}(\eta_{1,2},p')\in [\pi /3,\pi/2]$; this is 
the last claim (the angle can be bigger than $\pi /2$ if we have
exchanged $\eta_{i,j}$ and $\xi_{i,j}$ before; we still need to verify
that $p,p'$ have the other desired properties).

\medbreak

If $p'\in
\text{int}(F_{1,2}\cap F_{2,3})$, then it is immediate that
$(p',\eta_{2},\eta_{3},\eta_{4})$ is a tripod (because
$\overline{p'\eta_{1}}\cap \overline{p'\eta_{3}}\supsetneq \{p' \}$).

If $p'\in \partial
(F_{1,2}\cap F_{2,3})$, we still have (since $p'\in F_{2,3}$ by the
above and $p'\in F_{2,4}$ by definition)
\[
\angle_{p'}(\eta_{2},\eta_{4})=\pi =\angle_{p'}(\eta_{2},\eta_{3}). 
\]

If $p'$ is not tripodal for this triple, we would have to have 
$\angle_{p'}(\eta_{3},\eta_{4})\in \{0,\pi /3 \}$ (since the shift of
the triple is zero by assumption; see the proof of Proposition
\ref{existenceOfTripods}).  

However, $\angle_{p'}(\eta_{1},\eta_{4})=\pi$ by construction and 
$\angle_{p'}(\eta_{1},\eta_{3})\leq \pi /3$ since $p'$ cannot be
tripodal for $(\eta_{1},\eta_{2},\eta_{3})$.
Therefore, $\angle_{p'}(\eta_{3},\eta_{4})\geq 2\pi /3$, showing that
$p'$ is tripodal for the triple $(\eta_{2},\eta_{3},\eta_{4})$.

%

Similarly, we see $p\in F_{1,2}\cap F_{1,4}$ and that
$(p,\eta_{1},\eta_{3},\eta_{4})$ is a tripod.
%
%
\end{proof}

\bigbreak

The lemma above shows in particular that $F_{1,2}\cap F_{2,3}\cap
F_{3,4}\cap F_{1,4}\not =\emptyset$. Let us examine this set in more
detail, and give some more interpretation to the results from the
previous lemma:

\begin{lemma}\label{flatsIntersecting}
In the situation as in the previous lemma, we have 
\[
F_{1,2}\cap F_{3,4}=F_{1,4}\cap F_{2,3}=F_{1,2}\cap F_{2,3}\cap
F_{3,4}\cap F_{1,4}.
\]
\end{lemma}
\begin{proof}
Let us introduce a set $C$, drawn in Figure
\ref{flatsIntersectingFigure}: 
The left vertical boundary is 
\[
s_{1}:=l_{1,2,3}\cap l_{1,3,4}\subset
F_{1,2}\cap F_{2,3}\cap F_{3,4}\cap F_{1,4},
\]
which we have
just shown to be non-empty. 
Observe that every point in $s_{1}$ is tripodal for
$(\eta_{1},\eta_{2},\eta_{3})$ and for
$(\eta_{1},\eta_{3},\eta_{4})$. Every \emph{interior} point $x$ of
$s_{1}$ satisfies
$\overrightarrow{x\eta_{2}}=\overrightarrow{x\eta_{4}}$. 

\begin{figure}\def\color[#1]#2{}\def\mddefault{}
\def\updefault{}

\hfill 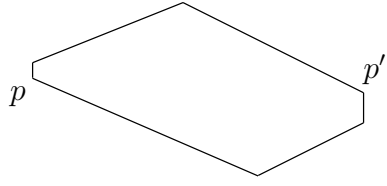\hfill ~
\caption{$F_{1,2}\cap F_{3,4}=F_{1,4}\cap F_{2,3}$}
\label{flatsIntersectingFigure}
\end{figure}

Similar properties hold for the vertical boundary on the right, which is
defined as 
\[
\emptyset \not =s_{2}:=l_{1,2,4}\cap l_{2,3,4}\subset F_{1,2}\cap F_{2,3}\cap
F_{3,4}\cap F_{1,4}. 
\]
Now we set $C$ to be the smallest convex polygon in $F_{1,2}$ 
containing $s_{1}\cup s_{2}$ and
such
that all the boundary segments are singular. (Observe that $C$ may
degenerate to a segment.)

%

By definition (and the two inclusions above), the set $C$ is a subset of
both $F_{1,2}\cap  F_{3,4}$ and of $F_{1,4}\cap F_{2,3}$.

Let us explain the relations to the previous lemma:
There, we have found that if $b_{1,2}(s_{1})\cap
b_{1,2}(s_{2})\not =\emptyset$, then there is an isometrically
embedded tree in $X$ with the given 4 asymptotic endpoints. 
If this is not the case, then we have made the assumption 
that $b_{1,2}(s_{1})>b_{1,2}(s_{2})$, and our choice of $p,p'$ was
such that $p$ is the lower endpoint of $s_{1}$ and $p'$ is the upper
endpoint of $s_{2}$. 

\medbreak
To finish the proof of our current lemma, we want to show that every boundary
segment of $C$ lies in the boundary of both $F_{1,2}\cap F_{3,4}$ and
$F_{1,4}\cap F_{2,3}$. 


This is immediate for the vertical segments $s_{1}$ and $s_{2}$. 

Observe that $b_{1}(s_{2})>b_{1}(s_{1})$ by the assumptions of Lemma
\ref{structureLemma}. 
Therefore, there are non-degenerate angular segments bounding $C$.
The argument for the angular boundary components are similar to each
other, let us give one in detail:

Let $\bar{p}$ be the upper endpoint of $s_{2}$, and
consider the segment $s:=\overline{\bar{p}\mu_{1}}\cap
C=\overline{\bar{p}\mu_{3}}\cap C$. 
Let us assume that $s$ is non-degenerate (i.e.~$s\not =\{\bar{p} \}$
and let $x$ be an interior point 
of~$s$. 

Note that $\bar{p}$ is the upper endpoint of $l_{1,2,4}$ or of
$l_{2,3,4}$.\footnote{In fact, the following argument shows that 
``$s$ non-degenerate implies that $\bar{p}$ is the upper endpoint of
both segments; see the remark below the lemma.} 
If $\bar{p}$ is the upper endpoint 
of $l_{1,2,4}$,  then
$s$ lies in the boundary of $F_{1,2}\cap F_{1,4}$; in particular,
we have 
\[
 \Sigma_{x}(F_{1,2}) \cap \Sigma_{x}(F_{2,3})
\ni\overrightarrow{x\mu_{2}}\not 
=\overrightarrow{x\mu_{4}}\in \Sigma_{x}(F_{3,4}) \cap \Sigma_{x}(F_{1,4}).
\footnote{Observe that we have
$\overrightarrow{x\bar{p}}=\overrightarrow{x\nu_{2}}=
\overrightarrow{x\nu_{4}}$; hence, we have 
$\angle_{x}(\eta_{2},\eta_{4})=\pi /3$.\label{tripodFootnote}}
\]
Otherwise, $\bar{p}$ is the upper endpoint of $l_{2,3,4}$, so
$s$ lies in the boundary of $F_{2,3}\cap F_{3,4}$; then 
the equation above holds as well.

The equation above shows that $x$ (and hence all of $s$) lies in the
boundary of both sets, 
$F_{1,2}\cap F_{3,4}$ and $F_{1,4}\cap F_{2,3}$.
%
Similar arguments hold for the other segments bounding $C$.
%
%
%
%
%
%
%
\end{proof}

\begin{remark}
Let us examine what the last argument shows about $\bar{p}$ (using the
notation of the previous lemma): Since $x\in C\subset F_{2,3}\cap F_{3,4}$,
the last footnote shows that there cannot be a tripod
$(\eta_{2},\eta_{3},\eta_{4})$ at $b_{2,3}$-level higher than
$b_{2,3}(\bar{p})$. Hence, $\bar{p}$ is the upper endpoint of
$l_{2,3,4}$. Similarly, $\bar{p}$ is the upper endpoint of
$l_{1,2,4}$. 

This shows that $l_{1,2,4}=l_{2,3,4}$ if the two angular boundary segments of
$C$ starting from the upper and lower endpoints of $s_{2}$ have
different slope.

The same statement holds for $l_{1,2,3}$ and $l_{1,3,4}$. 

Therefore, if $C$ has four angular boundary segments, then we have
$l_{1,2,4}=l_{2,3,4}$ and $l_{1,2,3}=l_{1,3,4}$.
\end{remark}

\begin{remark}
Let us also give a description of $C$ in terms of Busemann functions:
Normalize such that  
\[
b_{1}(p)=b_{2}(p)=b_{3}(p)=b_{4}(p)=0.
\]
Then we have for all $i\in \{1,2,3,4 \}$ (taking the indices modulo 4) 
\[
(b_{i}+b_{i+1})|_{F_{i,i+1}}=const = b_{i}(p)+b_{i+1}(p)=0.
\]
Now $x\in C$ if and only if $x\in F_{1,2}\cap F_{2,3}\cap F_{3,4}\cap
F_{1,4}$, as we have shown above. 
Let $f:=b_{1}+b_{2}+b_{3}+b_{4}$.
Then it follows that $x\in C$ implies 
\[
f(x)=
\frac{1}{2}((b_{1}(x)+b_{2}(x))+(b_{2}(x) +b_{3}(x)) +
(b_{3}(x)+b_{4}(x)) + (b_{1}(x)+b_{4}(x))) = 0.
\]
Since $x\in F_{i,i+1}$ if and only if $b_{i}(x)+b_{i+1}(x)=0$, and
$x\not \in F_{i,i+1}$ implies $b_{i}(x)+b_{i+1}(x)>0$, we have 
\[
x\in C \Leftrightarrow f(x)=0,
\]
so $C$ is the set of minima of $f$.
\end{remark}

\begin{lemma}\label{flatsNotIntersecting}
In the situation as in the previous lemmas, we have 
\[
F_{1,3}\cap
F_{2,4}=\emptyset.
\]
\end{lemma}
\begin{proof}
Choose points $p,p'$ as in Lemma \ref{structureLemma}, and normalize
(as above) such 
that $b_{1}(p)=b_{2}(p)=b_{3}(p)=b_{4}(p)=0$. 
Note that on $F_{i,j}$, we have $b_{i}+b_{j}=const$.

Let $x\in
\overline{p'\eta_{1}}$. Oberserve $\pi_{F_{2,4}}(x)=p'$. This implies
$b_{1}|_{F_{2,4}}\geq b_{1}(p')>b_{1}(p)=0$. 

Arguing similarly for $x\in \overline{p'\eta_{3}}$, we find
$b_{3}|_{F_{2,4}}\geq b_{3}(p') = 
-b_{2}(p')=b_{1}(p')>0$. Hence, $b_{1}+b_{3}|_{F_{2,4}}>0$, implying
the claim (since $b_{1}+b_{3}|_{F_{1,3}}\equiv 0$). 
\end{proof}

Let us phrase a version of Lemma \ref{structureLemma} 
which is valid for every S-set of cardinality 4:

\begin{proposition}\label{fourSingularPoints}
Let $\{\xi_{i}\st  i\in \{1,2,3,4 \} \}\subset \partial_{T}X$ be an S-set
of cardinality four. 
Set $\eta_{1}:=\xi_{1}$.
Then there are points $p,p'\in X$ and a numbering
$\{\eta_{2},\eta_{3},\eta_{4} \} = \{\xi_{2},\xi_{3},\xi_{4} \}$,
such that we have tripods 
\[
(p,\eta_{1},\eta_{2},\eta_{3}), (p,\eta_{1},\eta_{3},\eta_{4}),
(p',\eta_{1},\eta_{2},\eta_{4}), \text{ and }
(p',\eta_{2},\eta_{3},\eta_{4}).
\]
In particular, we have
\[
\overline{pp'}\subset F_{1,2}\cap F_{2,3}\cap F_{3,4}\cap F_{1,4}.
\]
\end{proposition}

\begin{proof}
If possible, we choose the identification
$\{\eta_{2},\eta_{3},\eta_{4} \} = \{\xi_{2},\xi_{3},\xi_{4} \}$ such that
\begin{align}\label{b1-eqn}
b_{1}(p_{1,2,3})\not =b_{1}(p_{1,2,4}).
\end{align}

Let us first assume that this is possible: then by exchanging
$\eta_{3},\eta_{4}$ if 
necessary, we may assume that $b_{1}(p_{1,2,3}) < b_{1}(p_{1,2,4})$.
Now Lemma \ref{structureLemma} applies (and finishes the proof).

We still need to consider the case that a choice as in \eqref{b1-eqn}
is not possible:
So pick an arbitrary identification $\{\eta_{2},\eta_{3},\eta_{4} \} =
\{\xi_{2},\xi_{3},\xi_{4} \}$, and
 assume that $b_{1}(p_{1,j,k})$ is independent of $j,k$.

We claim that in this case, there exists a 4-pod. 
By our assumptions, we have $b_{1}(l_{1,2,3})=b_{1}(l_{1,2,4}) =
b_{1}(l_{1,3,4})$.  If there is a point 
\[
p\in l_{1,2,3}\cap l_{1,2,4}\cap l_{1,3,4},
\]
then $p$ is the singular point of a 4-pod: This only means that $p\in
l_{2,3,4}$ as well, which follows immediately from the other three
inclusions. If this is the case (i.e.~if a 4-pod exists), then we set
$p'=p$, and we are done.

If the three sets above have pairwise non-empty intersection, then
they share a point.
Hence, we may assume that $l_{1,2,3}\cap l_{1,2,4}=\emptyset$. 

%
%
%
%
If this were the case,
the shift of
$(\eta_{2},\eta_{3},\eta_{4})$ cannot be 0.
The argument for this is the same as the one showing $p'\in F_{2,3}$
in the proof of Lemma \ref{structureLemma} (one can produce the tip
$p''$ of $F_{2,3}\cap F_{2,4}$, which should be a tripodal point, but
one can show that it cannot be).
%
%
%
\end{proof}

Let us summarize what we have achieved in this section:

Given a 4-point S-set $A$, there is either a 4-pod in $X$, or there is a
$2+2$ partition $A_{1}=\{a_{1},a'_{1} \},A_{2}=\{a_{2},a'_{2} \}$ of
$A$, such that 
\[
s_{1}:=l_{a_{1},a'_{1},a_{2}}\cap
l_{a_{1},a'_{1},a'_{2}}\not =\emptyset, \text{ and }
s_{2}:=l_{a_{1},a_{2},a'_{2}}\cap l_{a_{1}',a_{2},a'_{2}}\not
=\emptyset. 
\]
 In this case, the sets $s_{1}$ and $s_{2}$ can be joined
to each other ``almost horizontally'' (this is the statement about the
angle in Lemma \ref{structureLemma}).

\subsection{S-sets and trees}\label{treeSection}
Let $A:=\{\eta_{i}\st i\in I \}$ be an S-set. 
%
 
Let us examine the set $\F:=\bigcup_{i,j\in I}F_{i,j}$. 
We are going to construct a ``vertical'' quotient of $\F$ which is a
metric tree. 

Let $x\in F_{i,j}$, and consider some point $\eta_{k}\in A$.
Define 
\[
B_{k,i,j}(x):=b_{k}(\pi_{T_{i,j,k}}(x))=\min \Kr{b_{k}(\{y\in F_{i,j}\st
b_{i}(y)=b_{i}(x) \})}.
\]
\begin{remark}\label{defRemark}
\begin{enumerate}
\item If $k=i$ or $k=j$, the right-most definition still makes sense (and
$B_{k,i,j}(x)=b_{k}(x)$).
\item Note that the value of $B_{k,i,j}(x)$ does not depend on the choice of
$p_{i,j,k}$. 
\item If $y\in F_{i,j}\cap \{b_{i}=b_{i}(x) \}$, then
$B_{k,i,j}(y)=B_{k,i,j}(x)$. We will say that such a $y$
``represents'' $x$.
Using our convention for drawing flats $F_{i,j}$, this means that
``vertical'' lines all represent one point (in the space $\T$ which is
defined below).
\item Assume that $b_{i}(x)\leq b_{i}(p_{i,j,k})$. Then for any $y\in
F_{i,j}\cap F_{i,k}\cap \{b_{i}=b_{i}(x) \}$, we have
$B_{k,i,j}(x)=b_{k}(y)$. \label{usefulPoint}
\end{enumerate}
\end{remark}

We will see that the definition of $B_{k,i,j}(x)$ depends on
$x,\eta_{k}$ only (Lemma 
\ref{welldefined}), so we can define $B_{k}(x)$ (for $x\in \F$).

Now for every $k\in I$ the definition
\[
D_{k}(x,y):=|B_{k}(x)-B_{k}(y)|
\]
defines a pseudometric on $\F$; indeed, the triangle inequality follows
immediately from the inequality for real numbers.

We will see below that $D_{k}(x,y)\leq d(x,y)$.
Hence, the following is also a pseudometric on $\F$:
\[
D(x,y):=\sup_{k\in I}D_{k}(x,y)
\]
Consider the metric space $(\T,D):=(\F/\{D=0 \},D)$. 
In this section, we prove:

\begin{theorem}\label{tree}
$(\T,D)$ is a metric tree.
\end{theorem}

We start with some lemmas:
\begin{lemma}\label{stripsIntersecting}
Let $i_{0},i_{1},j_{0},j_{1}$ be four distinct elements of $I$. Then
\[
F_{i_{0},i_{1}}\cap F_{j_{0},j_{1}}\subset F_{i_{0},j_{0}}\cup
F_{i_{0},j_{1}}. 
\]
\end{lemma}
\begin{proof}
The claim is trivial if the intersection is empty.
Otherwise, it is an immediate consequence of Lemma
\ref{flatsIntersecting}. 
%
%
%
\end{proof}

\begin{lemma}
If $x\in F_{i,j}\cap F_{i,j'}$, then 
\[
B_{k,i,j}(x)=B_{k,i,j'}(x).
\]
\end{lemma}
\begin{proof}
If $k=i$, the claim is trivial. 

If $k=j$ (or analogously $k=j'$), we
have $B_{k,i,j}(x)=b_{j}(x) = B_{j,i,j'}(x)$ by Remark
\ref{defRemark}.\ref{usefulPoint}.  


So we may assume that $i,j,j',k$ are all distinct;
we consider the situation discussed in Lemma \ref{structureLemma}, and
assume that $\{i,j,j',k \}=\{1,2,3,4 \}$.  
We may assume $k=4$, and need to examine two cases: $i=1$ and $i=2$
(since $i=3$ is equivalent to $i=1$). 

If $i=2$, we have $x\in F_{1,2}\cap F_{2,3}$, and $b_{4}(\pi_{T_{1,2,4}}(x))
=b_{4}(\pi_{T_{2,3,4}}(x))$ follows: If 
$\pi_{T_{1,2,4}}(x)\in \overline{p'\eta_{2}}$,
the two projections are equal, and the claim follows. If not, 
we have $b_{2}(x)\in [b_{2}(p'),0]$, and we can represent $x$ in
$F_{1,2}\cap F_{2,3}\cap F_{3,4}\cap F_{1,4}$ (see Lemma
\ref{flatsIntersecting}). Now the claim follows from Remark
\ref{defRemark}.\ref{usefulPoint}. 

For  $i=1$, we have $x\in F_{1,2}\cap F_{1,3}$; by Lemma
\ref{structureLemma}, we have $\overline{p\eta_{1}}\subset F_{1,2}\cap
F_{1,3}\cap F_{1,4}$. Clearly, we can represent $x$ in this ray, and
the claim follows again from Remark \ref{defRemark}.\ref{usefulPoint}.
\end{proof}

\begin{lemma}\label{welldefined}
If $x\in F_{i,j}\cap F_{i',j'}$, then 
\[
B_{k,i,j}(x)=B_{k,i',j'}(x)=:B_{k}(x).
\]
\end{lemma}
\begin{proof}
If $\{i,j \}\cap \{i',j' \}\not =\emptyset$, then the claim follows
from the previous 
lemma. 
Otherwise, we may assume $x\in F_{i,j}\cap F_{i,j'}$ by 
Lemma \ref{stripsIntersecting}. Hence, we can apply  the previous
lemma twice:  
\[
B_{k,i,j}(x)=B_{k,i,j'}(x)=B_{k,i',j'}(x).\tag*{\qed}
\]
\renewcommand{\qed}{}
\end{proof}

We have shown that for every $k\in I$ the definition
\[
D_{k}(x,y):=|B_{k}(x)-B_{k}(y)|
\]
makes sense; so it is indeed 
 a pseudometric on $\F$ as claimed above.

Hence, the following is also a pseudometric on $\F$ (possibly with
value $\infty$):
\[
D(x,y):=\sup_{k\in I}d_{k}(x,y).
\]
Let $[x]$ denote the equivalence class of $x\in \F$. 
Recall that points $x,y\in F_{i,j}$ with $b_{i}(x)=b_{i}(y)$ satisfy
$[x]=[y]$. 

\begin{lemma}\label{essentialLemma}
Given points $x\in F_{i,j}, y\in F_{i',j'}$, there exist points
$x',y'\in F_{i'',j''}$ such that $\{i'',j'' \}\subset \{i,j,i',j' \}$ and
$[x]=[x'],[y]=[y']$, and
\[
D(x,y)= d(x',y').
\]
\end{lemma}
\begin{proof}
If $|\{i,j,i',j' \}|\leq 3$, we can project $x,y$ to a tripod or
line. In particular, we can represent $x,y$ by points $x',y'$ on a
line in a flat $F_{i'',j''}$.

If $|\{i,j,i',j' \}|=4$,
let us consider only the corresponding boundary points. We may
enumerate these as in Proposition 
\ref{fourSingularPoints}. 
Then we can represent $x$ and $y$ (uniquely) by points $x'',y''$ in 
$\overline{p\eta_{1}}\cup \overline{p\eta_{3}}\cup \overline{pp'}\cup
\overline{p'\eta_{2}}\cup \overline{p'\eta_{4}}$. Every two points in
this set lie in  a common flat, so let $x'',y''\in
F_{i'',j''}$. 

Choose 
\[
x'\in F_{i'',j''}\cap \{b_{i''}=b_{i''}(x'') \},\text{ and }
y'\in F_{i'',j''}\cap \{b_{i''}=b_{i''}(y'') \},
\]
such that
\[
d(x',y')=|b_{i''}(y'')-b_{i''}(x'')|.
\]
Now for all $k\in I$, we have $D_{k}(x,y)=D_{k}(x',y')\leq
d(x',y')$ (because projection to $T_{i'',j'',k}$ is
1-Lipschitz). Furthermore, we have
$D_{i''}(x,y)=D_{i''}(x',y')=d(x',y')$.  
\end{proof}

Hence, we have a metric space $(\T,D):=(\F/\{D=0 \},D)$. We claim that
$\T$ is a metric tree.

\begin{lemma}\label{fourPointTree}
If the cardinality of $A$ is 4, then $\T$ is a metric tree.
\end{lemma}
\begin{proof}
This is almost immediate from the previous lemma:
The discussion there
shows that  $\T$ has the topological structure of the set 
$\overline{p\eta_{1}}\cup \overline{p\eta_{3}}\cup \overline{pp'}\cup
\overline{p'\eta_{2}}\cup \overline{p'\eta_{4}}$ (assuming that the
elements of $A$ are named such that Proposition 
\ref{fourSingularPoints} and Lemma \ref{structureLemma} apply).
Now $D$ is almost the length metric on this graph: we just have to
shorten 
$\overline{pp'}$ to have length $b_{1}(p')-b_{1}(p)$. 
%
%
\end{proof}

%


\begin{proof}[Proof of Theorem \ref{tree}]
We put together the pieces collected above:
\begin{itemize}
\item For two points $x,y\in \F$, we can find a
flat $F_{i,j}$ and points $x',y'\in F_{i,j}$, such that
$[x]=[x'],[y]=[y']$, and $d(x',y')=D([x],[y])$ (Lemma
\ref{essentialLemma}). 

Then the segment $\overline{x'y'}$ represents a
geodesic $\overline{[x][y]}$ (of unit speed).
\item  From
Lemma \ref{fourPointTree}, we conclude that $\T$ has extendible
geodesics, and 
\item  Since for every $z\in \F$, the geodesics between $x',y',z$ (of
the form introduced 
above)  lie in a tree (again by Lemma \ref{fourPointTree}), every
triangle  in $\T$ is
degenerate. 
\item This implies that geodesic segments are unique, and that $\T$ is
0-hyper\-bolic. 
\end{itemize}
So $\T$ is indeed a tree.
\end{proof}

Let $\pi:\F\rightarrow \T$ be the projection, and observe that the
asymptotic endpoints of $\T$ correspond to the points $\eta_{i}$. We let
$\hat{\eta}_{i}$ denote the point of $\partial_{T}(\T)$ corresponding
to $\eta_{i}$. Then
$B_{i}(x)=b_{\hat{\eta}_{i}}([x])$.\footnote{Abusing notation, we will
sometimes also write $B_{i}([x]):=b_{\hat{\eta}_{i}}([x])$.}

\begin{lemma}\label{allFlatsIntersect}
Assume that $A$ is a good S-set, let $[x]\in \T$
, and let $\T_{[x]}:=\{(i,j)\st
[x]\in \pi (F_{i,j})=\overline{\hat{\eta}_{i}\hat{\eta}_{j}} \}$. 
Set 
\[
C_{[x]}:=\bigcap_{(i,j)\in \T_{[x]}}F_{i,j}.
\]
Then $C_{[x]}$ is non-empty, closed and convex,  and $[x]\in \pi (C_{[x]})$.
\end{lemma}
\begin{proof}
Let $(i_{0},j_{0})\in \T_{[x]}$. 
%

Consider $J:=\{j\st (i_{0},j)\in \T_{[x]} \}$.

For every $j,j'\in J$, we have 
\[
s_{j,j'}:=\sup\{b_{i_{0}}(F_{i_{0},j}\cap
F_{i_{0},j'} )\} =
b_{i_{0}}(p_{i_{0},j,j'})=B_{i_{0}}(p_{i_{0},j,j'})\geq B_{i_{0}}([x]). 
\]
The last inequality is due to the fact that both $(i_{0},j)$ and
$(i_{0},j')$ are in $\T_{[x]}$, hence $[x]\in
\overline{[p_{i_{0},j,j'}]\hat{\eta}_{i_{0}}}$. 

\begin{figure}\def\color[#1]#2{}\def\mddefault{}
\def\updefault{}

\hfill 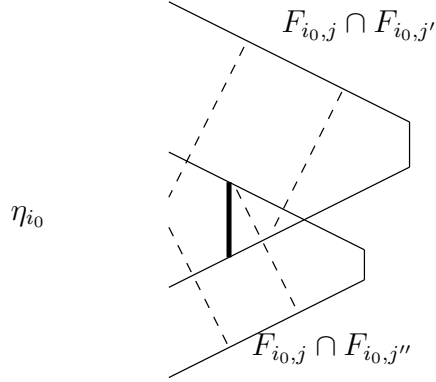\hfill\hfill  ~
\caption{The situation from Lemma \ref{allFlatsIntersect}. The fat line
is $l_{\{j,j',j'' \}}$. The tip of the inner sector is
$p_{i_{0},j',j''}$.} 
\label{ljjFigure}
\end{figure}

By induction, one shows:

For every finite $U\subset J$, the set 
$\{B_{i_{0}}=B_{i_{0}}([x])\}\cap \bigcap_{j\in U}F_{i_{0},j}$ 
is a non-empty geodesic segment $l_{U}$. Otherwise, we could find
$j,j'\in J$ such that $s_{j,j'}<B_{i_{0}}([x])$ (by
the argument for $S$ in the proof of Lemma \ref{structureLemma}). See Figure
\ref{ljjFigure}. 

If we can find $j,j'\in J$ such that $l_{\{j,j' \}}$ is compact, we
use this compactness to conclude that 
$l_{J}\not =\emptyset$: The sets $l_{\{j,j' \}}\backslash
l_{\{j,j',j'' \}}$ form an open cover of $l_{\{j,j' \}}$, so finitely
many $j''$ suffice, in contradiction to the above.

If such a choice of $j,j'$ is not possible, then the
assumption that $A$ is good implies that $l_{J}$ is a ray or a
geodesic line.

Similarly, let $J':=\{i\st (i,j_{0})\in \T_{[x]} \}$, and obtain
$l_{J'}':= \{B_{j_{0}}=B_{j_{0}}([x])\}\cap \bigcap_{i\in
J'}F_{i,j_{0}}$. 

If $l_{J}\cap l_{J'}'$ were empty, we could find $j\in J,i\in J'$ such
that $\{i_{0},j_{0},i,j \}$  contradict Proposition
\ref{fourSingularPoints}. 

Now $l_{J}\cap l_{J'}'\subset C_{[x]}$ by Lemma \ref{flatsIntersecting} (in
fact, $l_{J}\cap l_{J'}'=C_{[x]}\cap \{B_{i_{0}}=B_{i_{0}}([x])\}$).
%
%
%
\end{proof}

\begin{remark}\label{intersectingRemark}
If $I$ is finite or $X$ is discrete,  we can describe $C_{[x]}$ in detail as
follows: \newline
We find 
$j',j''\in J$ such 
that  we have $F_{i_{0},j'}\cap F_{i_{0},j''}\cap
\{B_{i_{0}}=B_{i_{0}}([x])\} = 
l_{J}$. 

Then we can similarly find $i',i''\in J'$ 
such that 
\[
l_{J}\cap l_{J'}'=
F_{i',j'}\cap F_{i'',j''}\cap \{B_{j'}=B_{j'}([x])\}
\]
(see also Lemma
\ref{flatsIntersecting} and Figure \ref{flatsIntersectingFigure}). To
cover a ``vertical cut-off'', we may have 
to introduce third indices $i''',j'''$ in $J',J$ resp., such that 
\[
C_{[x]}=F_{i',j'}\cap F_{i'',j''}\cap F_{i''',j'''}.
\]

In the general case, we can find sequences
$i_{n}',i_{n}'',i_{n}''',j_{n}',j_{n}'',j_{n}'''$ such that 
\[
F_{i_{n}',j_{n}'}\cap F_{i_{n}'',j_{n}''}\cap F_{i_{n}''',j_{n}'''}
\]
is a descending sequence with $C_{[x]}$ as its limit.
\end{remark}

\begin{remark}
Although the tree $\T$ is not isometrically embedded in $X$, the lemma
above shows that we can almost embed $\T$, and that intersection of
vertical lines in $\F$ is an equivalence relation. One may think of
the almost-embedded $\T$ in terms of sets as in Lemma
\ref{flatsIntersecting}. 
\end{remark}

%
%
%
%

The following lemma follows immediately from the definition:
\begin{lemma}\label{segmentConvex}
Let $[x],[y]\in \T$, and $C_{[x]},C_{[y]}$ from the lemma above. If
$[x],[y]$ and $\overline{[x][y]}$ are regular, then $C_{[x]}=C_{[y]}$.
\qed 
\end{lemma}
%
%
%
%
%
%
%
\bigbreak

At one point, we will need the following technical observation:

\begin{remark}\label{angleRemark}
Let $x\in C_{[x]}$, and 
$y\in C_{[y]}$. Assume that $x$ minimizes 
\[
b_{i',j'}|_{C_{[x]}\cap
\{b_{i'}=b_{i'}(x) \}}
\]
for some $(i',j')\in \T_{[x]}\cap
\T_{[y]}$. Then 
\[
\angle_{x}(y,\xi_{i',j'})\leq 2\pi /3. 
\] 
\emph{Reason:}
We may assume that $b_{j'}(y)\leq b_{j'}(x)$. Let $l$ be the line
joining $\nu_{j'}$ to $\mu_{i'}$ passing through $y$. Let $x'$ be the
point in $l$ satisfying $b_{j'}(x')=b_{j'}(x)$. 
Then it is easy to see that $b_{i',j'}(x')\geq b_{i',j'}(x)$.
\end{remark}

\section{Thickening tripods}
\label{thickeningSection}

Let $(p,\eta_{1},\eta_{2},\eta_{3})$ be a tripod in $X$.  We want to
find convex rank 
1-sets 
containing the tripod, other than a tubular neighborhood.

The results from this section are not used in the proof of Theorem
\ref{convBuild};
however, the techniques we introduce here are important for
the proof (and will be generalized later on).
Moreover, we get a feeling for the kinds of sets we use to build our
convex set later on.

We normalize the Busemann functions  to satisfy $b_{i,j}(p)=0$ for all
$i,j$. Let us agree to 
view the indices modulo 3.
Note that for the singular vertices $\eta_{i,j}$, the $\alpha$ used in
section \ref{busemannSection} is $\pi /3$.

Let us list some useful properties of the lower endpoint of
$l_{1,2,3}$:

\begin{lemma}\label{endpointLemma}
Assume that $p$ is the lower endpoint of $l_{1,2,3}$. Then (we will
list only one version, but permuting the indices 
leaves the statement intact, of course):
\begin{enumerate}
\item All the $\overrightarrow{p\eta_{i,j}}, \text{ }i,j\in \{1,2,3 \}$
 are distinct.
\item The $\overrightarrow{p\eta_{i,j}}$ span a flat in
$\Sigma_{p}X$. The singular directions of this flat are in the
directions of the $\eta_{i,j}$ and the $\nu_{i}$. In particular,
\label{flatProperty}
$\overline{p\nu_{1}}\cup \overline{p\eta_{2,3}}$ is a geodesic in $X$.
\item Let $x\in X\backslash \{p \}$. Then there exist $i',j'$ such that
$\angle_{p}(x, \eta_{i',j'})\geq 2\pi /3$. It follows that
$b_{i',j'}(x)\geq b_{i,j}(x)$ for all $i,j$.
\label{angleProperty}
\item If $b_{i,j}(x)\leq D$ for all $i,j$, then $d(x,p)\leq 2D$.
\label{distanceProperty}
\item If $b_{1,2}(x)>\max(b_{1,3}(x),b_{2,3}(x))$, then
$\angle_{p}(x,\nu_{3})<\pi /3$.
\label{angleProperty2}
\item If $b_{1,2}(x)>\max(b_{1,3}(x),b_{2,3}(x))$, then
$\angle_{q}(x,\nu_{3})<\pi /3$ for all $q\in l_{1,2,3}$.
\label{angleProperty3}
\item If $\angle_{p}(x,\eta_{1})\leq \pi /2$, we may distinguish two
cases:\label{twoCases}
\begin{enumerate}
\item $\angle_{p}(x,\nu_{1})\leq \pi /3$, which implies
$b_{2,3}(x)\geq b_{i,j}(x)$ for all $i,j\in \{1,2,3 \}$.  
\item otherwise $b_{1,2}(x)=b_{1,3}(x)>d(x,p)/2$, and
$\overrightarrow{x\eta_{1,2}}=\overrightarrow{x\eta_{1,3}}$. 
\end{enumerate}
\end{enumerate}
\end{lemma}
\begin{proof}
1: If two of the $\overrightarrow{p\eta_{i,j}}$ agree, then all three
have to be equal to each other; but then, $p$ is not the lower
endpoint of $l_{1,2,3}$.

Now 2 is clear.

3: Suppose $\angle_{p}(x,\eta_{1,j})<2\pi /3$ for both $j$. \newline
This is only possible if
$\angle_{p}(x,\nu_{1})<\pi /3$, and by~\ref{flatProperty}, we have
$\angle_{p}(x,\eta_{2,3})>2\pi /3$. The second part of the claim is
clear, since $b_{i',j'}$ increases at maximal slope along
$\overline{px}$ (in the sense of section \ref{busemannSection}). 

4:
By \ref{angleProperty}, at least one of the $b_{i,j}$
increases at maximal slope (at least $1/2 = -\cos(2\pi /3)$) along
$\overline{px}$. 

5: It follows from property 3, that
%
\[
\angle_{p}(x,\eta_{1,2}) \geq 2\pi /3>
\max(\angle_{p}(x,\eta_{1,3}),\angle_{p}(x,\eta_{2,3})).
\]
The claim follows as in the proof of 3.

6: We may assume $q\not =p$. Observe that
$\angle_{q}(x,\eta_{1,2})=\angle_{q}(x,\eta_{1,3})$, so this angle is
less than $2\pi /3$
(otherwise, $b_{1,3}(x)=b_{1,2}(x)$).

If $\angle_{q}(x,\eta_{1,2})=\angle_{q}(x,\eta_{1,3})=0$, let $q'$ be
the first point along $\overline{qp}$ where
$\angle_{q'}(x,\eta_{1,3})\not =0$; if such a point does not exist,
set  $q'=p$. Then 
$\alpha :=\angle_{q'}(x,\eta_{1,3})=2\pi /3$, or $q'=p$ and $\alpha =0$
(since the type of 
$\overrightarrow{q''x}$ does not change along $\overline{qx}$).
If $\alpha$ is $2\pi /3$, this is a contradiction to the above. If
$\alpha =0$,  
then $\angle_{p}(x,\eta_{2,3})=2\pi /3$, 
a
contradiction again. 

So there is a direction $\nu \in
\Sigma_{q}(X)$ of the 
same type as $\overrightarrow{q\nu_{3}}$ with
$\angle_{q}(\nu,\eta_{1,2})=\pi /3$ and $\angle_{q}(\nu ,x)<\pi
/3$. 
If $\nu \not =\overrightarrow{q\nu_{3}}$, this is a contradiction to
 Lemma \ref{triangleContradiction}. Hence, we have $\nu
=\overrightarrow{q\nu_{3}}$. 

7: First observe that $B_{\pi /2}(\overrightarrow{p\eta_{1}})=B_{\pi
/3}(\overrightarrow{p\nu_{1}})\cup B_{\pi
/3}(\overrightarrow{p\mu_{1}})$. \newline 
If $\angle_{p}(x,\nu_{1})\leq \pi /3$, case (a) follows immediately
(via \ref{flatProperty} \& \ref{angleProperty}).\newline
If this is not the case, then
$\alpha :=\angle_{p}(x,\eta_{1,2})=\angle_{p}(x,\eta_{1,3})>2\pi /3$,
implying the first part of the claim.
Since $\angle_{x}(\eta_{1,2},\eta_{1,3})\leq 2(\pi -\alpha )<2\pi /3$,
the second claim follows (because the two directions are singular and
of the same type).
\end{proof}

%
%
%
%

Pick $R>0$, and $D>R/2$.
We define convex sets as follows: Let $T_{1,2,3}:=\Conv
(p,\eta_{1},\eta_{2},\eta_{3})$ be the tripod, recall that we consider
the indices modulo 3, and let
\[
C_{i}:=B_{R}(T_{1,2,3})\cap \{b_{i,i+1}\leq D \}
\cap \{b_{i,i+2}\leq D \}.
\]

\begin{proposition}\label{firstThickening}
$C_{1}\cup C_{2}\cup C_{3}$ is convex.
\end{proposition}

We could prove this proposition directly; however, it can also
be derived from Proposition
\ref{joiningHoroballsComplicated}, so we omit a direct proof here.

To better understand the sets $C_{i}$, let us explain the
relation to the sets~$\tilde{C}_{i}$, which come to mind (more)
naturally; 
define 
\[
\tilde{C}_{i}:=B_{R}(\overline{p\eta_{i}})\cap \{b_{i,i+1}\leq D \}
\cap \{b_{i,i+2}\leq D \} = C_{i}\cap B_{R}(\overline{p\eta_{i}}).
\]

\begin{lemma}
If $p$ is the lower endpoint of $l_{1,2,3}$, we have
\[
\bigcup C_{i}=\bigcup \tilde{C}_{i}.
\]
\end{lemma}
\begin{proof}
The two sets in question are obviously equal on $B_{R}(p)$, but the
set on the left-hand side is potentially larger.
Consider a point $x\in X$ with $\pi_{T_{1,2,3}}(x)\in
\overline{p\eta_{1}}\backslash \{p \}$; note that points of
$B_{R}(\overline{p\eta_{1}})\backslash B_{R}(p)$ have this property.
We have $\angle_{p}(x,\eta_{1})<\pi /2$, so the cases from property
\ref{endpointLemma}.\ref{twoCases} apply. 

We see that 
%
in both cases, $x\in \bigcup C_{i}$ implies $x\in \tilde{C}_{1}$.
%
%
%
%
%
Since the conditions are symmetric, we are done.
\end{proof}


It will turn out that a convex rank 1-set as in Proposition
\ref{firstThickening} is not  quite good enough, 
so we need a 
 more
sophisticated approach:

In a first step, we show that we can do without tubular
neighborhoods, by imposing conditions on $b_{i,j}'$:

Normalize such that $b_{i,j}(p)=b_{i,j}'(p)=0$, let $D>0$ and $D'\in (D/2,2D)$. 
Consider the convex sets 
\begin{align*}
K_{i}:=\,&\{b_{i,i+1}\leq D \}\cap \{b_{i,i+2}\leq D \}\\
&\cap \{b_{i,i+1}'\leq D' \}\cap \{b_{i,i+2}'\leq D' \}
\end{align*}

\begin{proposition}\label{joiningHoroballsComplicated}
$\bigcup K_{i}$ is convex.
\end{proposition}

We start with an elementary observation:
\begin{lemma}\label{convexConditions}
\begin{align*}
\bigcup K_{i}\,&=\{x \st  \text{at least two }b_{i,j}(x)\leq D\}\\
&\quad ~\cap \{x\st \text{at least two }b_{i,j}'(x)\leq D' \}\\
&=:\tilde{K}_{1}\cap \tilde{K}_{2}
\end{align*}
\end{lemma}
\begin{proof}
If the claim is not true, then there is (\wlo) $x\in X$ with  
\[
b_{2,3}(x)>D\geq \max (b_{1,2}(x),b_{1,3}(x))
\text{ and }
b_{1,3}'(x)>D'\geq \max(b_{1,2}'(x),b_{2,3}'(x)).
\]
This implies that $l_{1,2,3}$ has a lower endpoint $p'$ and an upper
endpoint $p''$. We obtain
\[
\angle_{p'}(x,\nu_{1})<\pi /3
\text{ and }
\angle_{p''}(x,\mu_{2})<\pi /3
\]
by \ref{endpointLemma}.\ref{angleProperty2}.
This is a contradiction to Lemma \ref{triangleContradictionTwo}.
\end{proof}

\begin{proof}[Proof of Proposition \ref{joiningHoroballsComplicated}]
Note that $K_{1}\cap K_{2}=K_{2}\cap K_{3}=K_{3}\cap K_{1}$. 

We bring in the description $\bigcup K_{i}=\tilde{K}_{1}\cup
\tilde{K}_{2}$ from above:
We show convexity of $B_{\varepsilon}(x)\cap \tilde{K}_{i}$ for every
$x\in  K_{1}\cap K_{2}$, $i\in \{1,2 \}$, and an $\varepsilon >0$ that
we will 
construct in an instant. 
Via Lemma \ref{boundaryLemma} and the lemma above, this shows the claim.

It suffices to show convexity of $\tilde{K}_{1}$ near $x$,
since the proof is the same for $\tilde{K}_{2}$ (possibly with a 
different $\varepsilon$, but then Lemma \ref{boundaryLemma} applies to
the smaller one). 
\medbreak

\textbf{We construct $\varepsilon$:}\begin{itemize}
\item 
Pick $\varepsilon ,\hat{\alpha}$ such that in a 
Euclidean
triangle 
$\Delta (A,B,C)$ with $d(A,B)=2D$, $d(A,C)\in [2D-\varepsilon ,2D]$,
and $\angle_{A}(B,C)<\hat{\alpha}$, we have $d(B,C)<D\cdot
\sqrt{3}/2$; note that $\hat{\alpha}<\pi /3$.
\item 
We decrease $\varepsilon$ (if necessary), such that
\(
\varepsilon <\min (D\cdot
\sqrt{3}/2, (2D-D')/2).
\)
\item 
By decreasing $\hat{\alpha}$, we may assume that
$(2D-\varepsilon)\cdot \underbrace{\cos\hat{\alpha}}_{\approx  1}>D'$.
\item By
decreasing~$\varepsilon$ again, we can require
$(2D-\varepsilon 
)\cdot (\underbrace{-\cos(2\pi /3+\hat{\alpha })}_{> 1/2}) >D$. 
\end{itemize}

This is the $\varepsilon$ we work with.

Let $p'$ be the lower endpoint of $l_{1,2,3}$ (if $p'$ does not exist,
the claim for $\tilde{K}_{1}$ is trivial); then $b_{1,2}(p')\leq 0\leq
b'_{1,2}(p')$, and 
we set $R':=2\cdot (D-b_{1,2}(p'))\geq 2D$. Note that $K_{1}\cap
K_{2}\subset B_{R'}(p')$ (by
\ref{endpointLemma}.\ref{distanceProperty}).
Lemma~\ref{joiningHoroballs} shows convexity 
of $B_{R'}(p')\cap \tilde{K}_{1}$.

So it suffices to consider a point $x\in K_{1}\cap K_{2}$ with
$R'-\varepsilon\leq d(x,p')\leq R'$. Now for $i',j'$ from
\ref{endpointLemma}.\ref{angleProperty}, the construction of
$\varepsilon$ (last item) and $b_{i',j'}(x)\leq D$ imply
\[
\angle_{p'}(x,\eta_{i',j'})\in [2\pi /3,2\pi
/3+\hat{\alpha}].
\]
On the other hand, $b_{i,j}'(x)\leq D'$ implies that
$\angle_{p'}(x,\eta_{i,j})>\hat{\alpha}$ for all $i,j$ (by
construction of $\hat{\alpha}$)
,
so $\angle_{p}(x,\eta_{i,j})\in [2\pi
/3-\hat{\alpha},2\pi/3+\hat{\alpha}]$ for all $i,j$. 

This implies that there is a direction $\nu \in
\overline{\overrightarrow{p'\eta_{i',j'}}\overrightarrow{p'x}}$ such
that $\angle_{p'}(\nu ,\eta_{i,j})=2\pi /3$ for all $i,j$.

We can extend the flat half-strip $\Conv (x,p',\eta_{i',j'})$ to a flat
sector $F$ with tip $p'$, and inside this sector, we find a point $x'$
with $d(x',p')=R'$ and $\overrightarrow{p'x'}=\nu$. 
By construction of $\hat{\alpha}$, we have $d(x,x')<R'\cdot
\sqrt{3}/4$. Now 
Lemma~\ref{allAnglesMaximal} applies to $x'$. We have 
$B_{\varepsilon}(x)\subset B_{R'\cdot \sqrt{3}/2}(x')$, so this shows
the claim.
\end{proof}


This convex rank 1-set may have more asymptotic boundary points than
just the $\eta_{i}$.
We shrink it by putting in (large) tubular neighborhoods again:

Consider consistent (i.e.\ corresponding to each other under holonomy)
compact subsets $W_{i}$ 
 of $X_{\eta_{i}}$,
such that $[\overline{p\eta_{i}}]\in W_{i}$. Normalize such that
$b_{i,j}(W_{i})=[-S,S]=b_{i,j}'(W_{i})$. Let 
$S_{i,j}$ be the flat strip in $F_{i,j}$ ``spanned by'' $W_{i}$, i.e.
$S_{i,j}:=\Conv (W_{i},W_{j}) \subset F_{i,j}$. Let $R>10S$.

Let $\tilde{C}_{i}=S_{i,i+1}\cap S_{i,i+2}$, and let\footnote{instead
of $4S$, we could choose any value $D > 3S$
; the corresponding condition 
on $R$ would be $R>2(D+S)$.} 
\begin{align*}
C_{i}:=B_{R}(\tilde{C}_{i})&\cap \{b_{i,i+1}\leq 4S \}\cap
\{b_{i,i+2}\leq 4S \}\\ 
&\cap \{b_{i,i+1}'\leq 4S \}\cap \{b_{i,i+2}'\leq 4S \}\\
= B_{R}(\tilde{C}_{i})&\cap K_{i},
\end{align*}
where the $K_{i}$ are defined as before (with $D=4S+b_{1,2}(p),
D'=4S-b_{1,2}(p)\in [3S,5S]$, due to our new normalization).

\begin{proposition}\label{mostGeneralThickening}
$C:=\bigcup C_{i}$ is convex.
\end{proposition}

\begin{proof}
It suffices to show that $C':=C_{1}\cup C_{2}$ is convex. 

The last sentence above this proposition shows that
 Proposition \ref{joiningHoroballsComplicated} applies; hence
$C'\cap B_{R}(\tilde{C}_{1})\cap 
B_{R}(\tilde{C}_{2})$ is convex.

If some endpoint of $l_{1,2,3}$ lies in $S_{1,2}$, then $C_{1}\cap
C_{2}\subset B_{10S}(l\cap S_{1,2})$ (by
\ref{endpointLemma}.\ref{distanceProperty}), and since 
$R>10S$, we are done by Lemma \ref{boundaryLemma}. 

So we assume that no endpoint of $l_{1,2,3}$ lies in $S_{1,2}$.
Then the following lemma shows (in a precise way) that near
$C_{1}\backslash B_{R}(\tilde{C}_{2})$, the points in $C'$ lie in
$K_{1}$, and~$C'$ is convex in these points. Again, the claim
follows via Lemma \ref{boundaryLemma}. 
\end{proof}

\begin{lemma}\label{thickeningLemma}
Assume that no endpoint of $l_{1,2,3}$ lies in $S_{1,2}$. Then there exists an
$\varepsilon >0$ such that
if $x\in C_{2}$ and $y\in B_{2\varepsilon}(x)\cap C_{1}\backslash
B_{R}(\tilde{C}_{2})$, then $x\in K_{1}$, and
$\overline{xy}\subset C_{1}\cup C_{2}$.
\end{lemma}
\begin{proof}
Let us first construct the $\varepsilon$:\begin{itemize}
\item 
We pick $0<\hat{\alpha}<\pi /6$
such that $(R+10S)/2\cdot \underbrace{(-\cos(2\pi
/3-\hat{\alpha}))}_{<1/2}>5S$. 

\item 
Now pick $\varepsilon< (R-10S)/4$  such that in a Euclidean
triangle $\Delta (A,B,C)$ with 
\[
d(A,B)\in [R-4\varepsilon,R], d(A,C)\in [R-2\varepsilon
,R+2\varepsilon ], \text{ and }d(B,C)\in [0,2\varepsilon ],
\]
we have $\angle_{A}(B,C)<\hat{\alpha}$. 
\end{itemize}
This is the $\varepsilon$ (and $\hat{\alpha}$) we work with.

Now consider points $x,y$ as in the statement of the lemma.

The important step is the following observation: 

There exists a point $q\in \tilde{C}_{1}\cap
\tilde{C}_{2} = l\cap S_{1,2}=S_{1,2}\cap \{b_{1}=b_{1}(p) \}$ such that 
$\angle_{q}(x,\eta_{2})\geq \pi /2-\hat{\alpha}$.

\emph{Reason:}
If $\pi_{S_{1,2}}(x)\in \tilde{C}_{1}$, it is easy to pick $q\in
\tilde{C}_{1}\cap \tilde{C}_{2}$ suitably: set
$q:=\pi_{\tilde{C}_{2}}\circ \pi_{\tilde{C}_{1}}(x)$, and observe
$\angle_{q}(x,\eta_{2})>\pi /2$ (because
$\overline{\pi_{\tilde{C}_{1}}(x)} \cup \overline{q\eta_{2}}$ is a
geodesic ray). 

So assume that $\pi_{S_{1,2}}(x)\not \in \tilde{C}_{1}$. 
By definition, the point $y\in B_{2\varepsilon}(x)\cap C_{1}$ satisfies
$\pi_{S_{1,2}}(y)\in \tilde{C}_{1}$. Let $x'\in \overline{xy}$ such
that $q:=\pi_{S_{1,2}}(x')\in \tilde{C}_{1}\cap \tilde{C}_{2}$. 
Then $d(q,x')\in [R-4\varepsilon ,R]$, $d(q,x)\in [R-2\varepsilon
,R+2\varepsilon ]$, and $d(x,x')\leq 2\varepsilon$. 
Since $\angle_{q}(x',\eta_{2})\geq \pi /2$ by definition,
$q$ has the desired property by construction of $\varepsilon,
\hat{\alpha}$.  


\medbreak

Now assume that the claim is wrong. 
Then there is a point $x\in C_{2}$ as above with 
$\max(b_{1,2}(x),b_{2,3}(x))\leq 4S<b_{1,3}(x)$ (this is \wlo , maybe
we need to exchange $\eta_{2},\eta_{3}$ to get this inequality). 

As usual, we pick the lower endpoint $p'$ of $l_{1,2,3}$, for which we find
\[
\angle_{p'}(x,\nu_{2})<\pi /3
\]
by \ref{endpointLemma}.\ref{angleProperty2}. 
In particular, $\angle_{p'}(q,x)>\pi /3$.

Note that $\angle_{q}(p',x)>\pi /3+\hat{\alpha}$ (otherwise, we get
$b_{1,2}'(x)>4S$, because $b_{1,2}'(q)\geq -S$ and the second item in
the construction of $\varepsilon$). 
Now, if $\angle_{q}(x,\eta_{2})\in [\pi
/2-\hat{\alpha},\pi /2]$, 
we get $\angle_{q}(p',x)>2\pi/3-\hat{\alpha}$, 
implying $b_{1,2}(x)>4S$.

So $\angle_{q}(x,\eta_{2})>\pi /2$. 
Since $\angle_{p'}(q,x)>\pi /3$, we have $\angle_{q}(p',x)<2\pi
/3$. By the discussion above, there is a direction $\nu$ of the same
type as $\overrightarrow{q\nu_{1}}$, but neither
$\overrightarrow{q\xi_{1,2}}$ nor $\overrightarrow{q\nu_{2}}$, such
that  
\[
\angle_{q}(x,\nu)\leq \pi /3.
\]
But this is a
contradiction to Lemma \ref{triangleContradiction} (resp.~Remark
\ref{triangleContradictionRemark}). 

We have shown $x\in K_{1}$, so we have $x\in K_{1}\cap K_{2}$. If
$x\in B_{R}(\tilde{C}_{1})$, then $x\in C_{1}$ and the second claim is
immediate.  
If $x\not \in B_{R}(\tilde{C}_{1})$, then the argument from above,
applied to $y$, 
shows that $\{x,y \}\subset K_{1}\cap K_{2}$. Since
$B_{R}(\tilde{C}_{1}\cup \tilde{C}_{2})$ is convex, the second claim
follows. 
%
%
\end{proof}

\section{Existence of convex rank 1-sets}
\label{finalSection}
In this section, we prove Theorem \ref{convBuild}.
\subsection{Setting}
Let $A:=\{\eta_{i}\st i\in I \}$ be a good S-set.

For every triple $i,j,k\in A$, we pick a tripodal point
$p_{i,j,k}$. Let $S_{i,j,k}\in X_{\eta_{i}}$ be the strong asymptote
class at $\eta_{i}$ represented by
$\overline{p_{i,j,k}\eta_{i}}$. Since order of the indices does not
matter here, we can similarly define $S_{j,i,k}\in X_{\eta_{j}}$ and
so on.

Since all the shifts are 0, we can pick a particular $i_{0}\in I$, and
join all the strong asymptote classes $S_{i,j,k}$ to $\eta_{i_{0}}$,
where we obtain corresponding strong asymptote classes.

Let $K_{i_{0}}$ be the closed convex hull of all these strong asymptote
classes at~$\eta_{i_{0}}$. Since all the shifts are 0, we similarly
obtain isometric sets $K_{i}\subset X_{\eta_{i}}$ for all $i\in I$.

\begin{figure}\def\color[#1]#2{}\def\mddefault{}
\def\updefault{}

\hfill 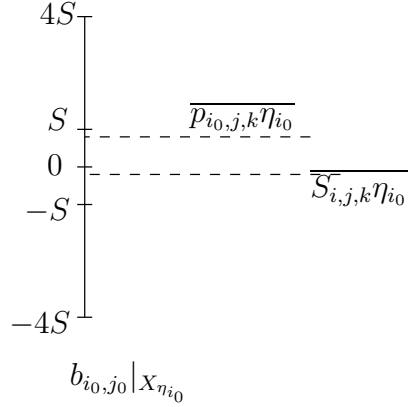\hfill ~
\caption{$X_{\eta_{i_{0}}}$ and the different kinds of elements of $K_{i_{0}}$.} 
\label{settingFigure}
\end{figure}

Because $A$ is good, we may assume that we have
chosen the $p_{i,j,k}$ such that the $K_{i}$ are compact.

We normalize the Busemann functions such that
\[
b_{i,j}(K_{i})=[-S,S]=b_{i,j}'(K_{i})
\]
(so we have $(b_{i,j}+b_{i,j}')|_{F_{i,j}}=0$.)


Recall the set $\F=\bigcup_{(i,j)\in I}F_{i,j}$, and its quotient tree
$\T$ from section \ref{treeSection}; as usual, we  let  
$\pi:\F\rightarrow T$ be the projection.


Also recalling the sets $\T_{[x]} = \{(i,j)\st [x]\in
\overline{\hat{\eta}_{i}\hat{\eta}_{j}} \}$, we 
set 
\[
\K_{[x]}:= 
\bigcap_{(i,j)\in \T_{[x]}} \{b_{i,j}\leq 4S
\}\cap \{b_{i,j}'\leq 4S  \}.
\]

In our choice of the limit $4S$, the important property is the
following: For every $p_{i,j,k}$, we have 
\[
4S-b_{i,j}(p_{i,j,k})\leq 5S<6S\leq 2(b_{i,j}(p_{i,j,k})-(-4S)).
\]
Of course, the same inequality holds for $b'_{i,j}$. These conditions
corresponds to the condition $D'\in (D/2,2D)$ in Proposition
\ref{joiningHoroballsComplicated}. Actually, one can extend both
results to the limit case where the inequality above is not strict;
however, this is not needed for the purpose of this paper.

\begin{lemma}\label{nonemptyLemma}
For every $[x]\in \T$, the set $\K_{[x]}$ is non-empty, closed, convex and
$[x]\in \pi (\K_{[x]})$.
\end{lemma}
\begin{proof}
Let $(i_{0},j_{0})\in \T_{[x]}$. 

We will use the notation of Lemma \ref{allFlatsIntersect}.

Clearly, it suffices to show that 
\[
\hat{C}_{[x]}:= \underbrace{C_{[x]}\cap
\{B_{i_{0}}=B_{i_{0}}([x])\}}_{=l_{J}\cap l_{J'}} \cap  
\{b_{i_{0},j_{0}}\in [-S,S] \} 
\not =\emptyset .
\]
If we have $b_{i_{0},j_{0}}(l_{J})<-S$, there is $j,j'\in J$ such that
$b_{i_{0},j_{0}}(l_{\{i_{0},j,j' \}})<-S$ (by Remark
\ref{intersectingRemark}), in contradiction to the 
construction of $K_{i_{0}}$.
Thus, we obtain $b_{i_{0},j_{0}}(l_{J})\cap [-S,S]\not
=\emptyset$ and $b_{i_{0},j_{0}}(l_{J'})\cap [-S,S]\not
=\emptyset$. Now the claim follows because $l_{J},l_{J'}$ are
intervals and have non-empty intersection.
%
%
\end{proof}

Let 
\[
\K:=\bigcup_{[x]\in \T}\K_{[x]}
\]

\begin{lemma}
$\K$ is connected.
\end{lemma}
\begin{proof}
Let $x\in \K_{[x]}, y\in \K_{[y]}$, and pick $(i',j')\in
\T_{[x]}\cap \T_{[y]}$.
We can join $x$ to $\hat{C}_{[x]}\subset S_{i',j'}$ and
$y$ to $\hat{C}_{[y]}\subset S_{i',j'}$.
By construction, we have $S_{i',j'}\subset \K$, so the claim follows.
\end{proof}

We are going to show that $\bar{\K}$ is 
 convex.
Since it is hard to show that $\bar{\K}$ is of rank 1, we introduce
tubular neighborhoods again:
Pick $R>10S$.
For $[x]\in \T$, let 
\begin{align*}
\tilde{C}_{[x]}\,&:=\K_{[x]}\cap B_{R}(\hat{C}_{[x]}),\\
\C\,&:=\bigcup_{[x]\in T}\tilde{C}_{[x]}.
\end{align*}

Exactly as for $\K$, we find that $\C$ is connected.
After showing that $\bar{\K}$ is convex, we also show that $\bar{\C}$
is convex.
Observe the analogon of moving from $\K$ to $\C$ 
and from
Proposition \ref{joiningHoroballsComplicated} to Proposition
\ref{mostGeneralThickening}.
For the new closed convex set $\bar{\C}$, it is easy to show that it
is of rank 1; this was obvious in both propositions mentioned above, 
because they were finite unions.
Thus, the proof of Theorem \ref{convBuild} is complete after these steps.

\subsection{The proof of Theorem \ref{convBuild}}


As a first step, we construct an $\varepsilon >0$, and show that
$\bar{\K}$ is 
$\delta$-locally convex for every 
$\delta <\varepsilon/2$. 

\begin{construction}\label{impEpsConstr}\begin{itemize}
\item
Pick $0<\hat{\alpha}<\pi /6$
such that 
\[
3S/\underbrace{(-\cos(2\pi /3+\hat{\alpha}))}_{>1/2}\geq
11S/2.
\]

\item 
By decreasing $\hat{\alpha}$ if necessary, we also require that
$11S/2\cdot \cos(\hat{\alpha}) > 5S$. 

\item 
Let $\varepsilon >0$ be such that in a
Euclidean triangle $\Delta (A,B,C)$ with 
\[
d(A,B)\geq 3S\text{ and }
d(B,C)\leq \varepsilon,
\]
 we have $\angle_{A}(B,C)<\hat{\alpha}/2$.
\end{itemize}
\end{construction}

We introduce some more notation for this section:
Consider points $[x_{0}]\not =[x_{1}]\in \T$. \newline
%
Pick $(i',j')\in \T_{[x_{0}]}\cap \T_{[x_{1}]}$
with
$B_{i'}([x_{0}])<B_{i'}([x_{1}])$.

Let $I_{0}:=\{i\in I\st (i,j')\in \T_{[x_{0}]} \}$. Analogously,
define $J_{0}:=\{j\in I\st (i',j)\in \T_{[x_{0}]} \}$ and
$I_{1},J_{1}$ (see Figure~\ref{treeFigure}). 
 Let $L:=J_{0}\cap I_{1}$.

Note that $\T_{[x_{0}]}\subset I_{0}\times J_{0}$, and 
$\T_{[x_{0}]}\cap \T_{[x_{1}]}=I_{0}\times J_{1}$. 

Set $\K_0:=\K_{[x_{0}]}, \K_{1}:=\K_{[x_{1}]}$. 

Let us start with a general lemma:
\begin{lemma}\label{convexityLemma}
Assume that $q\in \K_{[x_{0}]}\cap \K_{[x_{1}]}$ for
$[x_{0}],[x_{1}]\in \T$. Then $q\in \K_{[x]}$ for all 
$[x]\in \overline{[x_{0}][x_{1}]}$. 
\end{lemma}
\begin{proof}
Let $[x]\in \overline{[x_{0}][x_{1}]}$, and $(k,k')\in \T_{[x]}$ (see
Figure \ref{treeFigure}, with $[x]=[p_{x}]$). If 
one of $k,k'$ lies either in $I_{0}$ or in $J_{1}$, then
$b_{k,k'}(q)\leq 4S$ follows by assumption.
So we may assume $k,k'\in L$. 

It suffices to show $b_{k,k'}(q)\leq 4S$ for all $(k,k')\in L$, since
the 
claim for $b_{k,k'}'(q)$ follows analogously.

Assume that $b_{k,k'}(q)>4S$ for some $k,k'\in L$. 
Consider the lower endpoint $p$ of $l_{k,k',j'}$ and the lower
endpoint $p'$ of $l_{k,k',i'}$. By construction, we have
$B_{i'}(p')<B_{i'}(p)$. 

We have $b_{j',k}(q)\leq 4S$, $b_{j',k'}(q)\leq 4S$ and
$b_{k,k'}(q)>4S$. So
we obtain
$\angle_{p}(q,\nu_{j'})<\pi /3$  from
\ref{endpointLemma}.\ref{angleProperty2}. Similarly, we have 
$\angle_{p'}(q,\nu_{i'})<\pi /3$. 
By Lemma \ref{triangleContradiction} (and Remark
\ref{triangleContradictionRemark}),  
this is a contradiction.
%
\end{proof}


\begin{proposition}\label{segmentProposition}
Consider $x,y\in \K$ with 
$d(x,y)<\varepsilon$ (for the $\varepsilon$ from Construction
\ref{impEpsConstr}) and
 $x\in \K_0,y\in \K_{1}$
%
for some
$[x_{0}],[x_{1}]\in \T$.

Let $[q]\in \overline{[x_{0}][x_{1}]}$. Then 
\[
\K_{[q]}\cap \{x,y\} \not =\emptyset.
\]
\end{proposition}

\begin{proof}
We assume that $x\not \in \K_{[q]}$, and show that this implies $y\in
\K_{[q]}$.  
\Wlo, there is $(k_{x},k_{x}')\in \T_{[q]}$ with
$b_{k_{x},k_{x}'}(x)>4S$ (note that neither $k_{x}$ nor $k_{x}'$ lies
in $I_{0}$, because $x\in \K_{0}$). 
Pick $(i',j')\in \T_{[x_{0}]}\cap \T_{[x_{1}]}$ as above.

Let $p_{x}$ be the lower endpoint of $l_{i',k_{x},k_{x}'}$. 

We are
going to show $y\in \K_{[p_{x}]}$, which implies $y\in \K_{[q]}$ by
Lemma \ref{convexityLemma} (since $B_{i'}([p_{x}])\leq B_{i'}([q])$ by
construction). 

We have $b_{k_{x},i'}(x)\leq 4S$, $b_{k'_{x},i'}(x)\leq 4S$ and
$b_{k_{x},k'_{x}}(x)>4S$. So we have
\begin{align}\label{firstEqn}
\angle_{p_{x}}(x,\nu_{i'})<\pi /3\quad \text{by
}\ref{endpointLemma}.\ref{angleProperty2}. 
\end{align}
Let $p'$ be the lower endpoint of $C_{[p_{x}]}\cap
\{B_{i'}=B_{i'}(p_{x}) \}$. 
By (the proof of)
Lemma~\ref{nonemptyLemma}, $p'$ exists and satisfies
$b_{i,j}(p')\leq S$ for all $(i,j)\in \T_{[p_{x}]}$. 
In particular, we have $d(x,p')\geq 3S$.
This implies  
\begin{align}\label{secondEqn}
\angle_{p'}(x,y)< \hat{\alpha}/2 \text{ (by construction of $\varepsilon $)}.
\end{align}
We claim that we also have
\begin{align}\label{thirdEqn}
\angle_{p'}(x,\xi_{i',j'})<\pi -\hat{\alpha}.
\end{align}
Assume that this is not the case, and we have
$\angle_{p'}(x,\xi_{i',j'})\geq\pi -\hat{\alpha} (*)$.
Further, we have $\angle_{p'}(x,\eta_{k_{x},k'_{x}})<2\pi /3+\hat{\alpha}$.
Now
$b_{k_{x},k_{x}'}(x)>4S$ and
$b_{k_{x},k'_{x}}(p')=b_{k_{x},k'_{x}}(p')\leq S$ (the equality
follows from $p'\in C_{[p_{x}]}\subset F_{i',j'}\cap
F_{k_{x},k'_{x}}$); this implies 
$d(p',x)>11S/2$ (by construction of 
$\hat{\alpha}$). 

Taking $b_{i',j'}'(p')=-b_{i',j'}(p')\geq -S$ into account, $(*)$ and
$d(p',x)>11S/2$ imply $b_{i',j'}'(x)>4S$ (by construction of
$\hat{\alpha}$), in  
contradiction to $x\in \K_0$. Thus, \eqref{thirdEqn}  is proven.

\begin{figure}\def\color[#1]#2{}\def\mddefault{}
\def\updefault{}

\hfill 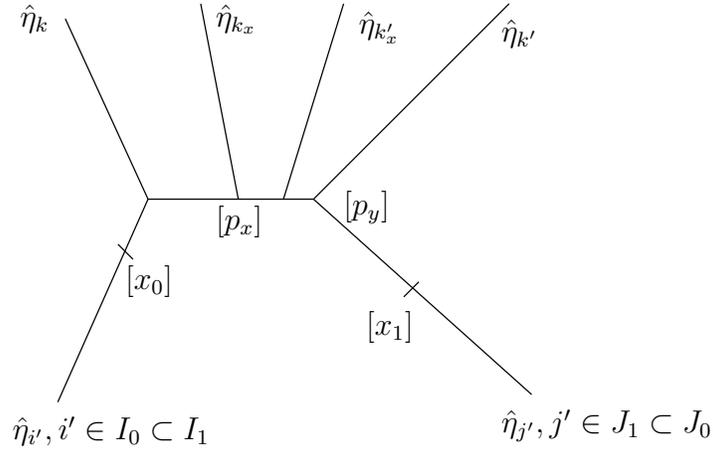\hfill ~
\caption{the relative position of the points in the tree $\T$. 
}
\label{treeFigure}
\end{figure}

Let us phrase the next steps as Lemmas:
\begin{lemma}\label{b-lemma}
We have $b_{k,k'}(y)\leq 4S$ for all $(k,k')\in
\T_{[p_{x}]}$. 
\end{lemma}
\begin{proof}
Assume that the claim is false, i.e.\ there are $(k,k')\in
\T_{[p_{x}]}$ with $b_{k,k'}(y)>4S$.
Observe that neither $k$ nor $k'$ lie in $J_{1}$, since $y\in \K_{1}$.
Let $p_{y}$ be the lower endpoint
of $l_{k,k',j'}$. 
We have 
\begin{align}
\angle_{p_{y}}(y,\nu_{j'})<\pi /3
\end{align}
by \ref{endpointLemma}.\ref{angleProperty2} (as in \eqref{firstEqn}).
As for \eqref{thirdEqn}, we obtain 
\[
\max
(\angle_{p_{y}}(y,\xi_{j',k}),\angle_{p_{y}}(y,\xi_{j',k'})) <
\pi -\hat{\alpha}. 
\]
Note that either $(j',k)$ or $(j',k')$ lie in $\T_{[p_{x}]}$ (so
$p'\in F_{j',k}$ or $p'\in F_{j',k'}$),
and that $B_{j'}(p_{y})\leq B_{j'}(p')=B_{j'}(p_{x})$. So Lemma
\ref{triangleContradiction} and Remark
\ref{triangleContradictionRemark} yield a contradiction (for
$p',p_{y},y$ and $\hat{\alpha}/2$).
\end{proof}

\begin{lemma}\label{b'-lemma}
We have $b_{k,k'}'(y)\leq 4S$ for all $(k,k')\in
\T_{[p_{x}]}$. 
\end{lemma}
\begin{proof}
Assume that this is not the case, i.e.\ there are $(k,k')\in
\T_{[p_{x}]}$ with $b_{k,k'}'(y)>4S$.
Observe that neither $k$ nor $k'$ lie in $J_{1}$, since $y\in \K_{1}$.
This time, let $p_{y}$ be the upper endpoint
of $l_{k,k',j'}$. 
We have 
\begin{align}
\angle_{p_{y}}(y,\mu_{j'})<\pi /3
\end{align}
by \ref{endpointLemma}.\ref{angleProperty2}.
As for \eqref{thirdEqn}, we obtain 
\[
\max
(\angle_{p_{y}}(y,\eta_{j',k}),\angle_{p_{y}}(y,\eta_{j',k'})) <
\pi -\hat{\alpha}. 
\]
Note that (at least) one of $(j',k)$ or $(j',k')$ lie in
$\T_{[p_{x}]}$, 
and that $B_{j'}(p_{y})\leq B_{j'}(p')$. 
We may assume that $(j',k)\in \T_{[p_{x}]}$ (by exchanging $k,k'$ if necessary). 

Since 
$p_{y}$ is the
\emph{upper} endpoint of $l_{k,k',j'}$, and $p'\in F_{k,k'}\cap
F_{j',k}\supset C_{[p_{x}]}$, we have
%
$\angle_{p_{y}}(p',\xi_{i',j'})\geq \pi /3$ (by Remark \ref{angleRemark}). 

So we have a contradiction to Lemma \ref{triangleContradictionTwo}. 
\end{proof}

This finishes the proof of Proposition \ref{segmentProposition}.
\end{proof}

\begin{proposition}\label{epsProp}
Let $x,y\in \K$ with
$d(x,y)<\varepsilon$ (for the $\varepsilon$ from Construction
\ref{impEpsConstr}).  
Then there exists $[q]\in \T$ such that $\overline{xy}\subset \K_{[q]}$.
\end{proposition}

\begin{proof}
As usual, let $x\in \K_{0},y\in \K_{1}$. 
By Lemma \ref{convexityLemma},
we know that the sets
\[
I_{x}:=\{[z]\in \overline{[x_{0}][x_{1}]} \st x\in \K_{[z]}
 \},
\]
and $I_{y}$ similarly for $y$, are intervals. By Proposition
\ref{segmentProposition}, $I_{x}\cup I_{y}$ 
covers $\overline{[x_{0}][x_{1}]}$. We want to show that $I_{x}\cap
I_{y}\not =\emptyset$.\footnote{If $X$ is discrete or $A$ is finite,
then the tree $\T$ is discrete. In this case, it is easy to see that both
$I_{x}$ and $I_{y}$ are open, so the claim follows.}

When we assume that this is not the case, then
we may assume that $I_{x}=\{[x_{0}]\}, I_{y}=
\overline{[x_{0}][x_{1}]} \backslash \{[x_{0}] \}$. 

Essentially, we want to show that $I_{x}$ is open; more specifically,
we will show that if $[b]\in \overline{[x_{0}][x_{1}]}$ is close
enough to $[x_{0}]$, then $x\in 
\K_{[b]}$. 

Pick $(i',j')\in \T_{[x_{0}]}\cap \T_{[x_{1}]}$ such that
$B_{j'}([x_{1}])< B_{j'}([x_{0}])$. 

Since $y\not \in \K_{0}$, there exist (\wlo) $(k,k')\in \T_{[x_{0}]}$
such that $b_{k,k'}(y)>4S$. 
We have $[p_{k,k',j'}]\in \overline{[x_{0}]\hat{\eta}_{j'}}$, so 
by (the proof of) Proposition
\ref{segmentProposition}, $x\in \K_{[p_{k,k',j'}]}$ holds, implying 
$[p_{k,k',j'}]=[x_{0}]$. 

Let $p$ be the lower endpoint of $\hat{C}_{[x_{0}]}$. 
By \ref{endpointLemma}.\ref{angleProperty3},
 we have 
\(
\angle_{p}(y,\nu_{j'})<\pi /3.
\)
By \cite[4.1.2]{kleinerLeeb}, there exists a point $a'\in
\overline{py}\backslash \{p \}$ such that $S':=\Conv (p,a',\xi_{i',j'})$
is a flat half-strip and $\overline{a'\xi_{i',j'}}\cap
\overline{p\mu_{j'}}\not =\emptyset$. 

Similarly, there exists a point $a''\in
\overline{py}\backslash \{p \}$ such that $S'':=\Conv (p,a'',\nu_{j'})$ 
is a flat half-strip.

Pick $a\in \text{int}(
\overline{pa'}\cap \overline{pa''})$, 
Then by construction, we have 
\[
\angle_{a}(y,\nu_{j'})<\pi /3.
\]

As for \eqref{thirdEqn}, we find $\angle_{p}(y,\xi_{i',j'})< \pi
-\hat{\alpha}$. 
%
Since $\angle_{a}(y,\xi_{i',j'})=\angle_{p}(y,\xi_{i',j'})$, the point
$a$ has the same  property (which we will need in order to apply Lemma
\ref{triangleContradiction}).

Now let
$\{b \}:=\overline{a\xi_{i',j'}}\cap \overline{p\mu_{j'}}$. 
This point exists by construction and lies in $F_{i',j'}$. So
$\K_{[b]}$ is defined.

Observe that $B_{j'}(b)<B_{j'}(p)= B_{j'}([x_{0}])$ by construction, so
$x\not \in \K_{[b]}$. 

\textbf{We claim that $x\not \in \K_{[b]}$ leads to a contradiction,
which finishes the proof.} 

\noindent
\textbf{Step 1:} $b_{k_{x},k_{x}'}(x)\leq
4S$ for all  $(k_{x},k_{x}')\in \T_{[b]}$. 

Assume that $b_{k_{x},k_{x}'}(x)>4S$ for some $(k_{x},k_{x}')\in \T_{[b]}$.
Let $p'$ be the lower endpoint of $\hat{C}_{[p_{k_{x},k_{x}',i'}]}
\subset F_{i',j'}$. By construction, we have 
\[
B_{j'}(b)\leq
B_{j'}(p')<B_{j'}(p)
\]
(the last inequality follows from
$(k_{x},k_{x}')\not \in \T_{[x_{0}]})$, and by
\ref{endpointLemma}.\ref{angleProperty3}, we have  
\[
\angle_{p'}(x,\nu_{i'})<\pi /3.
\]
We claim that Lemma
\ref{triangleContradiction} leads to a contradiction (for $a,p',y$ and
$\hat{\alpha}/2$; as in the proof
of Lemma \ref{b-lemma}). This is clear if $a\in F_{i',j'}$.\newline
 If $a\not
\in F_{i',j'}$, then $\angle_{p}(a,\xi_{i',j'})>2\pi /3$, but 
$\angle_{p}(p',\xi_{i',j'})\leq 2\pi /3$ (if $b_{i',j'}(p)=-S$, this
is trivial, because $p'\in S_{i',j'}$; otherwise, it follows from
\ref{angleRemark}).  
Therefore, $p'\in S'$ (because $B_{j'}(p')\geq B_{j'}(a)$), so we can
apply Lemma 
\ref{triangleContradiction} as claimed.

\noindent
\textbf{Step 2:}
$b_{k_{x},k_{x}'}'(x)\leq
4S$ for all  $(k_{x},k_{x}')\in \T_{[b]}$. 

Assume that $b_{k_{x},k_{x}'}'(x)>4S$ for some $(k_{x},k_{x}')\in \T_{[b]}$.
This time, let $p'$ be the upper endpoint of
$\hat{C}_{[p_{k_{x},k_{x}',i'}]} 
\subset F_{i',j'}$. As before, we have $B_{j'}(b)\leq
B_{j'}(p') < B_{j'}(p)$, and by
\ref{endpointLemma}.\ref{angleProperty3}, we have  
\[
\angle_{p'}(x,\mu_{i'})<\pi /3.
\]
We have $\angle_{p}(p',\xi_{i',j'})\leq 2\pi /3$ as above.
If $\angle_{p'}(a, \xi_{i',j'})\geq \pi /3$, we can apply Lemma
\ref{triangleContradictionTwo} (as in the proof of Lemma
\ref{b'-lemma}). 
Otherwise, we have
 $a\in F_{i',j'}$ and $\overrightarrow{ay}\in
\overline{\overrightarrow{a\mu_{j'}}\overrightarrow{a\nu_{j'}}}\subset
\Sigma_{a}(X)$.  
In this case, Lemma \ref{triangleContradiction} applies as in Step 1
(after exchanging the $\nu_{j}$ with the $\mu_{j}$).

Together, steps 1 and 2 show that $x\in \K_{[b]}$, the desired
contradiction. 
\end{proof}

Proposition \ref{epsProp} says: Whenever we consider $x,y\in \K$ with
$d(x,y)<\varepsilon$, then $\overline{xy}\subset \K$. This property is
inherited by the closure $\bar{\K}$.
This implies that $\bar{\K}$ is
$\delta$-locally convex for every $\delta <\varepsilon/2$. From Proposition
\ref{convexityLocalProperty}, we obtain:

\begin{theorem}
$\bar{\K}$ is convex.\qed 
\end{theorem}

%
%
%
%
%
%

It is hard to decide whether $\bar{K}$ is of rank 1. Hence, we bring in
additional conditions again: Pick $R>10S$.
For $[x]\in \T$, recall the set $\hat{C}_{[x]}$ from Lemma
\ref{nonemptyLemma}, and let 
\begin{align*}
\tilde{C}_{[x]}\,&:=\K_{[x]}\cap B_{R}(\hat{C}_{[x]}),\\
\C\,&:=\bigcup_{[x]\in T}\tilde{C}_{[x]}.
\end{align*}
As for $\K$, we find that $\C$ is connected. 

We want to show that $\bar{\C}$ is convex, by the same tools as for
$\bar{\K}$: 

\begin{proposition}\label{epsPropC}
There exists $\varepsilon >0$ such that
for $x,y\in \C$ with
$d(x,y)<\varepsilon$, 
we have  $\overline{xy}\subset \C$.
\end{proposition}

\begin{proof}
We pick $\varepsilon ,\hat{\alpha}$ such that they satisfy the
conditions from the proof of Lemma
\ref{thickeningLemma} as well as those from
Construction \ref{impEpsConstr}; this is possible, because in both
constructions, we first impose conditions on $\hat{\alpha}$, and
afterwards, we require $\varepsilon >0$ to be small enough.

Assume that $x\in \tilde{C}_{[x_{0}]},y\in \tilde{C}_{[x_{1}]}$. We know from
Proposition \ref{epsProp} that there is $[q]\in
\overline{[x_{0}][x_{1}]}$ with $\overline{xy }\subset \K_{[q]}$. 
If $\{x,y \}\subset \tilde{C}_{[q]}$, there is nothing to show.

Assume that $x\not \in \tilde{C}_{[q]}$: We know that $x\in \K_{[z]}$
for all $[z]\in \overline{[x_{0}][q]}$ by Lemma \ref{convexityLemma}. 
Hence, we have 
\begin{align}\label{newEqn}
\{[z] \in \overline{[x_{0}][q]}\st x\in
\tilde{C}_{[z]}\} = \{[z] \in \overline{[x_{0}][q]}\st x\in
B_{R}(\hat{C}_{[z]})\}.
\end{align}
Lemmas \ref{nonemptyLemma} and 
\ref{allFlatsIntersect}  imply that $\hat{C}_{[z]}$
varies continuously along~$\overline{[x_{0}][q]}$. 

Therefore (by pushing $[x_{0}]$ towards $[q]$ as far as possible), we
may assume $x\not 
\in 
\tilde{C}_{[z]}$ for every $[z]\in
\overline{[x_{0}][q]}\backslash \{[x_{0}] \}$, and
$d(x,\hat{C}_{[x_{0}]})=R$ $(*)$. 

Let $(i',j')\in \T_{[x_{0}]}\cap \T_{[x_{1}]}$ such that
$B_{j'}([x_{1}])< B_{j'}([x_{0}])$ (as usual). For every singular $[z]\in
\overline{[x_{0}][q]}$, we have
\begin{align}\label{stripEqn}
b_{i',j'}(\hat{C}_{[z]})=[-S,S],
\end{align}
since otherwise, $d(x,\hat{C}_{[z]})\leq
10S<R$ (by \ref{endpointLemma}.\ref{distanceProperty}), implying $x\in
\tilde{C}_{[z]}$.  

Similarly, we may assume $y\not \in \tilde{C}_{[z]}$ for every $[z]\in
\overline{[q][x_{1}]}\backslash \{ [x_{1}]\}$, and we
 get \eqref{stripEqn}  for every
singular $[z]\in \overline{[q][x_{1}]}$. 

Recalling from Figure \ref{flatsIntersectingFigure} what the sets
$\F_{[z]}$ look like, we may conclude that
%
$\bigcup_{[z]\in \overline{[x_{0}][x_{1}]}}\hat{C}_{[z]}$ is
convex (a convex subset of the strip $S_{i',j'}=\Conv (K_{i'},K_{j'})$; not
necessarily a rectangle, if $[x_{0}]$ and/or $[x_{1}]$ are not
singular), and so is 
\[
\bigcup_{[z]\in
\overline{[x_{0}][x_{1}]}}B_{R}(\hat{C}_{[z]})
= B_{R}(\bigcup_{[z]\in \overline{[x_{0}][x_{1}]}}\hat{C}_{[z]})
.
\]

Along the lines of Lemma \ref{thickeningLemma}, we obtain $y\in
\K_{[x_{0}]}$ and similarly $x\in \K_{[x_{1}]}$ (see below).  
Then it is immediate (from Lemma \ref{convexityLemma} and convexity of
the metric) that $\overline{xy}\subset \C$.

Let us explain the argument for $y\in \K_{[x_{0}]}$:

Assume that $y\not \in \K_{[x_{0}]}$, so
\wlo, we have $b_{k,k'}(y)>4S$ for some $(k,k')\in \T_{[x_{0}]}$. 

Consider the lower endpoint $p'$ of $l_{k,k',j'}$. It satisfies
$\angle_{p'}(y,\nu_{j'})<\pi /3$ by \ref{endpointLemma}.\ref{angleProperty3}.


Let $x':=\pi_{\hat{C}_{[x_{0}]}}(x)$, and
$\{y' \}:=\{\pi_{\hat{C}_{[p']}}(x') \}=\overline{x'\eta_{j'}}\cap
\hat{C}_{[p']}$.

If
$\angle_{y'}(y,\eta_{j})>\pi /2$, we get a contradiction to either the
sum of angles in a triangle, or Lemma \ref{triangleContradiction}.

Observe that $d(y',x)\geq  R$ and 
$\angle_{y'}(x,\eta_{j'})\geq \pi /2$ (because
$\angle_{x'}(x,\eta_{j'})\geq \pi /2$ by~$(*)$ and \eqref{newEqn}). This implies
$d(y',y)\geq R-\varepsilon$ and 
$\angle_{y'}(y,\eta_{j})\geq \pi /2-\hat{\alpha}$.

Now we obtain a contradiction as in the proof of Lemma \ref{thickeningLemma}.
\end{proof}

Just as for $\bar{\K}$, we now obtain that $\bar{\C}$ is convex. We claim that
it is also of rank 1.

\begin{theorem}
$\bar{\C}$ is a convex
rank 1-subset of $X$. 
\end{theorem}
\begin{proof}
If $\partial_{T}\bar{\C}$ is not a 0-dimensional subbuilding, then there
exists (\wlo) a point $\xi_{i,j}\in \partial_{T}\bar{\C}$. 
In fact, by \cite{centers}, $\partial_{T}\bar{\C}$ is a subbuilding or
has a center. So either all $\eta_{i,j},\xi_{i,j}$ are in the
asymptotic boundary, or all $\xi_{i,j}$ agree (again \wlo; it could
also be the $\eta_{i,j}$ that agree).

So consider a point $x\in F_{i,j}$ with
$b_{i,j}(x)=  S+2R+3\varepsilon$ (for some $\varepsilon >0$). Let
$x':=\pi_{\hat{C}_{[x]}}(x)$. Then 
$d(x,x')\geq 2R+3\varepsilon$. 

To finish the proof, it suffices to lead the following assumption to a
contradiction: There 
exists $[x'']\in \T$ such that $x\in
B_{\varepsilon}(\tilde{C}_{[x'']})$. 

Assume the contrary, and set
$x'':=\pi_{\hat{C}_{[x'']}}(x)$. Obviously, $d(x',x'')\geq 
R + 2\varepsilon$. Pick $i'$ such that $(i',j)\in \T_{[x]}\cap
\T_{[x'']}$ 
(such an $i'$ exists, after exchanging $i,j$ if necessary). 

Since $x',x''\in S_{i',j}$, the inequality $2S<(R+2\varepsilon )/2$
implies that 
\[
\angle_{x'}(x'',\xi_{i,j}) = \angle_{x'}(x'',x)\geq \pi /3.
\]

Now triangle comparison yields 
$d(x'',x)\geq R+2\varepsilon$, the
desired contradiction.
\end{proof}

By definition and Lemma \ref{nonemptyLemma}, we have $S_{i,j}\subset
\bar{\C}$ for all $i,j\in 
I$. Hence, we have $A\subset \partial_{T}\bar{\C}$. We have just shown
that $\bar{\C}\subset B_{2R}(\bigcup_{i,j\in I}S_{i,j})$. Therefore, 
%
$\partial_{T}\bar{\C}$
is precisely the closure of $A\subset \partial_{\infty}X$ in the cone
topology. 
The proof of Theorem \ref{convBuild} is now finished.

 \end{document}